\documentclass{article}
\usepackage{amsmath}
\usepackage{graphicx,psfrag,epsf}
\usepackage{enumerate}
\usepackage{natbib}
\usepackage{tabularx} 
\usepackage{amssymb}
\usepackage{amsfonts}
\usepackage{arydshln}
\usepackage{eufrak}
\usepackage{color}
\usepackage{bigstrut}
\usepackage{multirow}
\usepackage{bm}
\usepackage {tikz, calc}
\usetikzlibrary {positioning}
\usepackage[a4paper]{geometry}
\usepackage[textwidth=8em,textsize=small]{todonotes}
\numberwithin{equation}{section}

\usepackage{natbib}

\newtheorem{thm}{Theorem}
\newtheorem{lem}{Lemma}

\newtheorem{cor}{Corollary}

\begin{document}

\begin{center}\Large{Construction of Blocked Factorial Designs to Estimate Main}

\Large{Effects and Selected Two-Factor Interactions}

\vspace{0.2in}

\Large{J.D.Godolphin}

\vspace{0.2in}

\large{Department of Mathematics, University of Surrey, UK, GU2 7XH.}

\end{center}
\vspace{0.2in}

\noindent{\bf Abstract}
Two-level factorial designs are widely used in industrial experiments. For processes involving \(n\) factors, the construction of designs 
comprising \(2^n\)  and \(2^{n-p}\) factorials, arranged in blocks of size \(2^q\) is investigated. The aim is to estimate all main effects and a selected subset of two-factor interactions. Designs constructed according to minimum 
aberration criteria are shown to not necessarily be the most appropriate 
designs in this situation. A design construction approach is proposed which 
exploits known results on proper vertex colourings  in graph theory. Examples are provided to illustrate the results and construction strategies. 
 Particular consideration is given to the special case of designs with blocks of size four and tables of designs are given for this block size.

\vspace{0.2in}


\section{Introduction}
\label{sec:intro}

In full and fractional factorial experiments, sources of systematic variation
can be controlled by grouping runs into blocks of \(2^q\) heterogeneous 
experimental units  enabling the estimation of effects with greater 
precision. For blocking to be effective, blocks should generally be relatively 
small. Mead {\it et al.}\,(2012, Chapter 15) suggest that for many 
experiments, blocks should contain no more than eight experimental units. 
However, the arrangement of a factorial experiment in small blocks puts 
significant constraints on its estimability capability, due to the large number 
of factorial effects confounded with blocks. Thus, the construction of such 
designs with desirable properties can be challenging. 

We use the notation  \(2^{n-p}\) factorial for a regular fractional factorial in \(n\) two-level factors and with \(2^{n-p}\) runs. A \(2^{n-p}\) factorial is 
completely determined by \(p\) independent defining words, which generate the treatment defining contrast sub-group. The smallest number of factors involved in a word in this sub-group is the 
{\it resolution} of the design. Estimability capability can vary considerably between
 \(2^{n-p}\) factorials with the same resolution. The capacity of a design to provide estimates can be quantified under the assumptions of the effect hierarchy principle, namely:
\begin{enumerate}
\item[H1:] For \(m_1<m_2\), \(m_1\)-factor interactions are more likely to be important than \(m_2\)-factor interactions.
\item[H2:] For given \(m\), all \(m\)-factor interactions are equally important. 
\end{enumerate}
See Wu and Hamada (2009, Chapter 4).  Under assumptions H1 and H2, the minimum aberration criterion, introduced by Fries and Hunter (1980), ranks 
\(2^{n-p}\) factorials by the distribution of word lengths in the treatment defining contrast sub-group.
Tables of designs according to this criterion are provided by Chen {\it et al.} (1993) and Block and Mee (2005).
Several papers, including Hedayat and Pesotan (1992), Wu and Chen (1992) and Wang (2007),  address the problem of selecting  \(2^{n-p}\) factorials in situations where assumption H2 does not hold in that estimates of main effects and selected two-factor interactions are required. 
Wu {\it et al.} (2012) give tables of resolution IV \(2^{n-p}\) factorials for up to 128 runs which cover all estimable configurations of two-factor 
interactions. 

When arranging \(2^n\) factorials in blocks, two assumptions are made:
\begin{enumerate}
\item[B1:]  Interactions between block factors and treatment factors are negligible.
\item[B2:] An interaction between two or more block factors has the same importance as the main effect of a block factor.
\end{enumerate}
The blocking framework for a \(2^n\) factorial arranged in blocks of size 
\(2^q\) is typically described by means of a blocking defining contrast 
sub-group. This has the same mathematical structure as the treatment defining contrast sub-group for a \(2^{n-(n-q)}\) factorial, and is determined by 
\(n-q\) independent defining words. Sun {\it et al.}\,(1997) apply the 
minimum aberration criterion to the block defining contrast sub-group and  
give optimal blocking schemes for \(n \le 8\). 

The problem of selecting blocked \(2^{n-p}\) factorials with good 
estimability capability is considerably more complex than that of selecting a 
\(2^{n-p}\) factorial or of blocking a \(2^{n}\) factorial, being complicated by the need to take account of both treatment and block defining contrast 
sub-groups. A substantial body of work focuses on construction of optimal 
blocked \(2^{n-p}\) factorials. Sun {\it et al.}\,(1997) and Mukerjee and Wu 
(1999) use admissibility criteria based on both word length distributions. 
Other work involves various methods of combining the word length 
distributions of the two defining contrast subgroups into a single 
distribution, and using the minimum aberration criterion on this combined 
distribution.
See for example  Sitter {\it et al.}\,(1997), Chen and Cheng (1999), Zhang 
and Park (2000), Cheng and Wu (2002), Xu and Lau (2005)  and Xu and 
Mee (2010).  Instead of employing a word length approach, Cheng and 
Mukerjee (2001) use a criterion appertaining to the alias pattern of the 
interactions.  Several of the papers referred to above include extensive tables of 
designs. In particular, Sun {\it et al.} (1997) constitutes a comprehensive 
design catalogue for \(n \le 9\) and up to 128 runs.

Prior knowledge about the process under investigation can result in a subset of the two-factor interactions being of particular  interest.
The aim of this work is to develop a strategy for the bespoke arrangement 
of \(2^n\) and \(2^{n-p}\) factorials in blocks of size \(2^q\) so that all 
main effects and selected two-factor interactions can be estimated. We make assumptions H1, B1 and B2, but digress from the  effect hierarchy 
principle in that some two-factor interactions are considered more important than others.  In many cases a blocked full or fractional factorial will not 
provide estimates of all main effects and all two-factor interactions. 
A design approach which prioritises estimation of main effects and has 
practitioner led focus on those two-factor interactions thought most likely to 
be relevant therefore seems  pragmatic. No two-factor interactions are  
assumed  negligible: the estimation of any of these in addition to the 
specified ones is advantageous. Interactions in three or more factors are assumed to be negligible. Henceforth, unless otherwise stated,  {\it 
interaction} will refer to a two-factor interaction.

Designs in the literature comprising \(2^n\) and \(2^{n-p}\) factorials in blocks of size \(2^q\) yielding estimates of all main effects and some  interactions  are generally chosen to maximise the number of estimable  interactions, under assumption H2. 
If a specific set of  interactions are of interest, it may not be possible to find a design in 
the literature which accommodates all main effects and all chosen 
interactions. This can arise for two reasons. First,  it might be impossible to 
construct a design providing the desired estimates. This depends on both the number and configuration of the selected interactions.
Second, it may be possible to construct a design which accommodates the 
required estimates, but if such a design does not rank highly under 
optimality criteria then it will not appear in design tables. Further, even if a 
selected set of interactions can be accommodated by a listed design, the 
problem of mapping the factors of an experiment onto the factors of a 
design so all chosen interactions are estimable is not necessarily trivial.

In this work, instead of using the block defining contrast sub-group, we identify the blocking structure of a design by the 
sub-group of treatment combinations comprising the principal block. With 
blocks of size \(2^q\), this requires \(q\) independent treatment combinations, rather than \(n-q\) or \(n-p-q\) independent defining words, which is advantageous 
when \(q\) is small compared to \(n-q\) or \(n-p-q\). Further, these independent treatment combinations correspond to a \(q \times n\) generator matrix which gives 
valuable and immediate  information on the estimability of main effects and  interactions.  

In \S2, fundamental concepts on generator matrices are introduced and a graph theory approach is developed for the construction of designs comprising \(2^n\) 
factorials in blocks of size \(2^q\). By representing the selected interactions graphically, known results in graph theory yield conditions on the interaction set which guarantee the 
existence of a design that provides the required estimates.  For interaction sets satisfying the conditions, the construction problem is expressed as a vertex 
colouring task in graph theory. A valid colouring leads to a {\it factor grouping} associated with a {\it profile set}, which in turn yields a generator matrix and 
hence a suitable design. Thus, rather than seeking an existing design to use for an experiment, a design is constructed which is tailored to the set of chosen interactions. Graph theory terminology  is consistent with Bondy and Murty (2008).
In \S3, the method is extended to encompass constructions of  blocked \(2^{n-p}\) factorials. 
In \S4 we focus on the special case of designs with blocks of size four. This constitutes an important design class since the blocks are small enough for blocking to be effective in most cases, but large enough to allow the 
estimation of all main effects and some two-factor interactions. 
Design templates are provided for \(2^{n-p}\) factorials in blocks of size four for up to twelve factors and 128 runs. 

\section{Blocked full factorial designs}
\label{full}
\subsection{Preliminaries}\label{prelim}

Treatment combinations are represented in two equivalent forms: in lower 
case letters and as  row vectors \((x_{1}, x_{2}, \cdots, x_{n})\) in the Galois Field of order 2, GF(2), with  \(x_{j}
=1\) if the \(j\)th factor is at high level, and \(x_{j}=0\) otherwise. The 
treatment combination with all factors at low level is denoted by \((1)\) or by 
the \(1 \times n\) vector with all terms zero.  Factors are labelled \(F_1, 
\cdots, F_n,\) or alphabetically \(A, B, \ldots\) in examples. Thus,  in a five 
factor process, the treatment combination with \(A,C,D\) high and \(B,E\) 
low is represented as \(acd\) and \((1,0,1,1,0)\).

A blocked \(2^{n}\) factorial has  the  treatment combinations arranged in 
\(2^{n-q}\) blocks of size \(2^q\). The blocking structure can be identified 
by \(q\) treatment combinations represented by \(q\) independent row vectors denoted \(\bm{x_i}=(x_{i1}, x_{i2}, \cdots, x_{in})\), for \(i=1,\ldots,q\). 
The principal block contains the \(2^q\) treatment combinations of the form 
\(\sum_{i=1}^q\alpha_i \bm{x_i}\), where \(\alpha_i \in \{0,1\}\) and addition is modulo 2.
These form an Abelian sub-group of the \(2^n\) treatment 
combinations and remaining blocks contain cosets of this sub-group. Thus the design is fully determined by the \(q \times n\) generator matrix of rank \(q\):
\begin{equation*}\label{X}
X=\left( \begin{array}{c}
\bm{x_1}\\
\vdots \\
\bm{x_q}\\
\end{array} \right)=
\left( \begin{array}{cccc}
x_{11} & x_{12} & \cdots & x_{1n} \\
\vdots & \vdots & & \vdots \\
x_{q1} & x_{q2} & \cdots & x_{qn} \end{array} \right).
\end{equation*}  

\noindent{\bf Example~1:}~{\it
Consider design \(\mathcal{D}1\), comprising a 
\(2^5\) factorial arranged in blocks of size four with generator matrix:
\begin{equation*}
X_1=\left( \begin{array}{ccccc}
1 &1 & 1 &0& 0 \\
1& 0&1&1 & 1\end{array} \right).
\end{equation*} 
The principal block of \(\mathcal{D}1\) contains \((1)\), \(abc\), \(acde\) and \(bde\).  The design is given below, with columns representing blocks.
\begin{equation}\label{D1}
 \begin{array}{cccccccc}
(1) &a & b &c& d &e&ad&ae\\
abc&bc &ac&ab&abcd&abce&bcd&bce\\
acde& cde& abcde & ade& ace & acd& ce & cd\\
bde & abde & de & bcde & be & bd & abe & abd\end{array}.
\end{equation}}

In a blocked \(2^n\) factorial, the \(2^n-1\) factorial effects separate into two sets: \(2^{n-q}-1\) effects which are confounded with blocks, and 
\(2^n-2^{n-q}\) estimable effects. Estimability of main effects and interactions can be deduced directly from the generator matrix.

\begin{thm}\label{thm1}
A factorial effect is estimable from a blocked \(2^{n}\) factorial {\it iff} at least one of the generating treatment combinations contains an odd number of effect terms.
\end{thm}
\noindent{\bf Proof:}~Consider any factorial effect. Treatment combinations with an odd number of terms in common with the effect have one sign in the effect contrast, and treatment combinations with an even number of terms in common with the effect have the opposite sign. 
Thus, for the effect of the principal block to sum out in the linear combination of observations relating to the factorial effect, \(2^{q-1}\) treatment combinations of the principal 
block must have an odd number of effect terms and \(2^{q-1}\) must have an even number. This occurs {\it iff} at least one of the generating treatment 
combinations contains an odd number of effect terms. By the method of construction, every block contains \(2^{q-1}\) treatment combinations of each sign in the 
effect contrast, enabling all block effects to be eliminated, {\it iff} this property is exhibited by the principal block. This establishes the result. 

\vspace{0.1cm}

Theorem \ref{thm1} leads to results on estimability of main effects and interactions in terms of properties of the design generator matrix. To state the results efficiently, we introduce the notation \(\mathcal{X}_u\) to denote the set of \(2^{u}-1\) non-zero \(u \times 1\) vectors over GF(2). For example:
\begin{equation*}
\mathcal{X}_2=\left\{\left( \begin{array}{c}1\\1\\ \end{array} \right)\!,\!
\left( \begin{array}{c}1\\0\\ \end{array} \right)\!,\!
\left( \begin{array}{c}0\\1\\ \end{array} \right)
\right\}.
\end{equation*}
From Theorem \ref{thm1}, a necessary and sufficient condition for a main effect, \(F_j\) say, to be estimable is that at least one of \(x_{ij}\) 
is non-zero, for \(i=1, \ldots,q\). Thus,
\begin{cor}\label{cor1}
All main effects are estimable from a blocked \(2^{n}\) factorial {\it iff} every column of \(X\)  is a column of \(\mathcal{X}_q\).
\end{cor}

In this work, only designs from which all main effects are estimable will be 
considered, that is, all blocked \(2^n\) factorials are required to satisfy 
Corollary \ref{cor1} and therefore to have generator matrices consisting of 
columns of \(\mathcal{X}_q\). Further, since generator matrices are full rank, every generator matrix must contain at least \(q\) linearly independent vectors of \(\mathcal{X}_q\). 
Now consider an interaction, \(F_jF_k\) say. By Theorem 
\ref{thm1}, this is estimable {\it iff} there exists some \(i \in \{1, \ldots,q\}\) for which exactly one of \(x_{ij}\) and \(x_{ik}\) is non-zero. This 
gives:
\begin{cor}\label{cor2}
An interaction is estimable from a blocked \(2^{n}\) factorial {\it iff} the corresponding two columns of \(X\) are different.
\end{cor}

\vspace{0.2cm}

\noindent{\bf Example~1 continued:}~{\it
By Corollary \ref{cor1}, all five main effects are estimable from \(\mathcal{D}1\). By Corollary \ref{cor2}, \(\mathcal{D}1\) provides estimates of eight of the 
ten interactions. The first and third columns of \(X_1\) are the same indicating that \(AC\) is not estimable. Similarly, \(DE\) is not estimable. These properties are confirmed by inspection of \eqref{D1}.}


\vspace{0.2cm}

Corollaries \ref{cor1} and \ref{cor2} give an indication of the limitations of 
the capacity of a blocked \(2^n\) factorial to estimate interactions. From 
Corollary \ref{cor1}, there are \(2^q-1\) different columns available for inclusion in \(X\). From Corollary \ref{cor2} an interaction is only estimable if the corresponding columns of \(X\) are different. Note that for a design in 
blocks of size two, that is with \(q=1\),  the implications of Corollaries \ref{cor1} and \ref{cor2} 
are that if all main effects are estimable then no two-factor interactions are 
estimable, since \(\mathcal{X}_1\) only contains one vector. In general, if 
\(n>2^q-1\) not all  interactions will be estimable. The number of factors in 
the design can be expressed as \(n=(2^q-1)v+w\), where 
\(v=[n/(2^q-1)]\), with \([.]\) denoting {\it the integer part of}, and 
\(w \in \{0,1,\ldots, 2^q-2\}\). 
\begin{thm}\label{int}
An upper bound for the number of estimable interactions from a blocked \(2^{n}\) factorial is:
\begin{equation}\label{bound}
\phi_{\max}=\binom{n}{2}-vw -(2^q-1) \binom{v}{2}.
\end{equation} 
\end{thm}
\noindent{\bf Proof:}~A blocked \(2^n\) factorial has generator matrix \(X\) 
with columns comprising members of \(\mathcal{X}_q\). From Corollary 
\ref{cor2}, the number of estimable interactions is the number of pairs of 
columns of \(X\) in which the columns are different. Expressing the number 
of columns of \(X\) as \(n=(2^q-1)v+w\), an upper bound for the number of pairs of distinct columns in 
\(X\)  is:
\begin{equation*}
\binom{n}{2}-w\binom{v+1}{2}-(2^q-1-w)\binom{v+1}{2}.
\end{equation*}
After manipulation of the Binomial coefficients involving \(v+1\), this expression gives the right hand side of \eqref{bound} and the result follows. 

\vspace{0.1cm}

A similar result is given in Godolphin (2019) in relation to factorial designs arranged in a rectangular array. Theorem \ref{int} prompts:

\vspace{-0.2cm}

\begin{cor}\label{cor3}
If \(n\le 2^q-1\) then all interactions in a blocked \(2^{n}\) factorial are estimable, provided that  \(X\) is formed from \(n\) distinct members of 
\(\mathcal{X}_q\)  and has rank \(q\), over GF(2).
If \(n>2^q-1\) the number of estimable interactions in a blocked \(2^{n}\) factorial  achieves \(\phi_{\max}\) iff \(X\) is formed from \(v +1\) copies of each of \(w\) columns 
of \(\mathcal{X}_q\) and of \(v\) copies of each of the remaining \(2^q-1-w\) columns.
\end{cor}
Corallary \ref{cor3} is illustrated by the following examples:

\vspace{0.1cm}

\noindent{\bf Example~2:}~{\it A \(2^6\) factorial is to be arranged in 
blocks of size eight. Since \(n=6 < 7=2^q-1\), all main effects and interactions are estimable from a design with generator matrix having columns comprising six distinct members of \(\mathcal{X}_3\). Note that any six members of 
\(\mathcal{X}_3\) yield a full rank generator matrix. Such a matrix is 
\begin{equation*}
X_2=\left( \begin{array}{cccccc}
1 &0 & 0 &1& 1&0 \\
0& 1&0&1 & 0 &1\\
0 & 0 & 1&0&1&1 \end{array} \right),
\end{equation*}
giving design, \(\mathcal{D}2\) say, with principal block containing the 
treatment combination subgroup \((1), ade,bdf,\) \(abef,cef,\)\,\(acdf,\)\,\(bcde\) and \(abc\). The seven remaining blocks contain cosets.}

\vspace{0.1cm}

\noindent{\bf Example~3:}~{\it For a \(2^6\) factorial in blocks of size four, \(\phi_{\max}=12\). This is achieved by \(\mathcal{D}3\) with 
a generator matrix formed from two copies of each column of \(\mathcal{X}_2\):
\begin{equation*}
X_3=\left( \begin{array}{cccccc}
1 &1 & 0 &0& 1&1 \\
0& 0&1&1 & 1 &1\end{array} \right).
\end{equation*} 
The principal block of \(\mathcal{D}3\) contains \((1),\,abef,\,cdef\) and \(abcd\), and the 15 remaining blocks contain cosets.   Design \(\mathcal{D}3\) provides estimates of all  interactions except \(AB,\,CD\) and \(EF\). }

\vspace{0.1cm}

There are situations in which a design providing estimates of fewer than \(\phi_{\max}\) interactions is preferable to one achieving the bound. 

\vspace{0.1cm}

\noindent{\bf Example~4:}~{\it  A \(2^6\) factorial in blocks of size four is required 
to assess all main effects and all  interactions involving one factor, \(A\) say.  
Use of \(\mathcal{D}3\) is unsatisfactory since, even with a reordering of 
the columns of \(X_3\), one of the interactions of interest will be 
inestimable. An alternative design, \(\mathcal{D}4\), is suggested, with 
generator matrix:
\begin{equation*}
X_4=\left( \begin{array}{cccccc}
1 &0 & 0 &0& 1&1 \\
0& 1&1&1 & 1 &1\end{array} \right).
\end{equation*}
Design \(\mathcal{D}4\) gives estimates of 11  interactions, including all 
those involving \(A\).}

\subsection{Graphical representation of required and estimable  interactions}\label{graph}

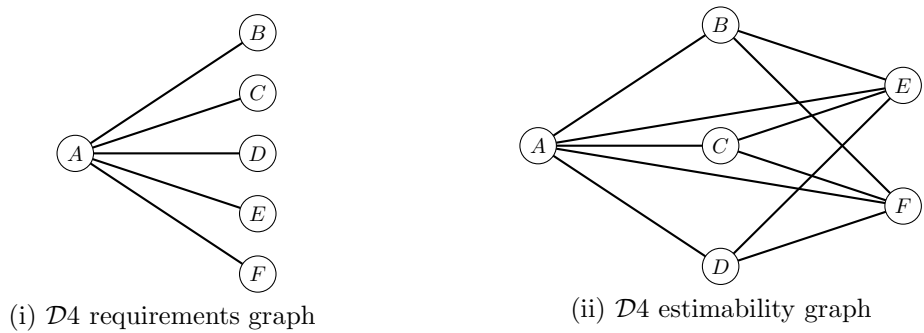
\begin{figure}
\caption{Graphs for  Example 4.}\label{eg4D4}
\begin{minipage}{.5\linewidth}
    \centering
\vspace{0.2cm}
\hspace{0cm} \begin{tikzpicture}[scale=0.8, transform shape] 
   \begin{scope} [vertex style/.style={draw,
                                       circle,
                                       minimum size=6mm,
                                       inner sep=0pt,
                                       outer sep=0pt}] 
\node[vertex style] (1) at (0,7) {$A$};
\node[vertex style] (2) at (3,9) {$B$};
\node[vertex style] (3) at (3,8) {$C$};
\node[vertex style] (4) at (3,7) {$D$};
\node[vertex style] (5) at (3,6) {$E$};
\node[vertex style] (6) at (3,5) {$F$};
    \end{scope}
 \begin{scope} 
\draw[thick] (1)--(2);
\draw[thick] (1)--(3);
\draw[thick] (1)--(4);
\draw[thick] (1)--(5);
\draw[thick] (1)--(6);
\end{scope}
\end{tikzpicture}

(i) \(\mathcal{D}4\) requirements graph
\end{minipage}%
\begin{minipage}{.5\linewidth}
\centering
\begin{tikzpicture}[scale=0.8, transform shape]
\begin{scope} [vertex style/.style={draw,
                                       circle,
                                       minimum size=6mm,
                                       inner sep=0pt,
                                       outer sep=0pt}] 
\node[vertex style] (1) at (0,7) {$A$};
\node[vertex style] (2) at (3,9) {$B$};
\node[vertex style] (3) at (3,7) {$C$};
\node[vertex style] (4) at (3,5) {$D$};
\node[vertex style] (5) at (6,8) {$E$};
\node[vertex style] (6) at (6,6) {$F$};
    \end{scope}
 \begin{scope}
\draw[thick] (1)--(2);
\draw[thick] (1)--(3); 
\draw[thick] (1)--(4);
\draw[thick] (1)--(5);
\draw[thick] (1)--(6);
\draw[thick] (2)--(5);
\draw[thick] (2)--(6);
\draw[thick] (3)--(5);
\draw[thick] (3)--(6);
\draw[thick] (4)--(5);
\draw[thick] (4)--(6);
\end{scope}
\end{tikzpicture}  

(ii) \(\mathcal{D}4\) estimability graph
\end{minipage}%
\end{figure}

When planning an experiment it is informative to represent interactions of interest by means of a labelled graph with vertices labelled to correspond to 
factors and with an edge between vertices \(F_i\) and \(F_j\) if interaction \(F_iF_j\) is of interest.  Hedayat and Pesotan (1992), Wu and Chen (1992) and Wu {\it et al.}\,(2012) use 
such graphs to aid  in the construction of designs to estimate main effects and selected interactions in experiments without blocking. Wu and Chen (1992) 
coin the phrase {\it requirements graph} which we adopt in this work. Requirements graphs have the advantages of giving a visual depiction of the way in which interactions of interest are interlinked and of 
making results in graph theory available as tools in the design planning process.  Figure \ref{eg4D4}(i) depicts the requirements graph for the problem of Example 4. The graph of  estimable interactions for a design will be termed the {\it estimability graph}. Figure \ref{eg4D4}(ii) gives the estimability 
graph for \(\mathcal{D}4\). This contains the five edges of the requirements graph and six additional 
edges, indicating that \(BE\),\,\(BF\),\,\(CE\),\,\(CF\),\,\(DE\) and \(DF\) are estimable as well as the selected  interactions. 
The estimability graph is a graphical representation of \(X_4\) since each edge coincides with a pair of generator matrix columns that are not equal. Note that the estimability graph has the requirements graph as a sub-graph.

\vspace{0.1cm}

The estimability graph for  \(\mathcal{D}4\), is a  {\it labelled complete 3-partite graph}.
In general, a  labelled complete \(r\)-partite graph on \(n\) vertices has the \(n\) labelled vertices assigned to \(r\) non-empty disjoint sets. There is an edge 
between every pair of vertices from different sets and there are no edges between vertices from the same set. In \(\mathcal{D}4\) the three vertex sets are \(\{B,C,D\}\), \(\{E,F\}\) and \(\{A\}\).
We denote the set of 
labelled complete \(r\)-partite graphs on \(n\) vertices, for \(r\) taking all values in \(\{q,q+1, \ldots, 2^q-1\}\),  as \(T(n,q)\). 
Every design comprising a \(2^n\) factorial in blocks of size \(2^q\) 
corresponds to a unique estimability graph in \(T(n,q)\). For example, the estimability graph for \(\mathcal{D}4\), given in Figure \ref{eg4D4}(ii), is a member of 
\(T(6,2)\). 

\begin{figure}
\caption{Graphs for Example 4 with vertex colouring.}\label{eg4D4a}
\begin{minipage}{.5\linewidth}
    \centering
\vspace{0.2cm}
\hspace{0cm} \begin{tikzpicture}[scale=0.8, transform shape] 
   \begin{scope} [vertex style/.style={draw,
                                       circle,
                                       minimum size=6mm,
                                       inner sep=0pt,
                                       outer sep=0pt}] 
                                       \node[vertex style] (1) at (0,7) {$A{\bm 3}$};
\node[vertex style] (2) at (3,9) {$B{\bm 1}$};
\node[vertex style] (3) at (3,8) {$C{\bm 1}$};
\node[vertex style] (4) at (3,7) {$D{\bm 1}$};
\node[vertex style] (5) at (3,6) {$E{\bm 2}$};
\node[vertex style] (6) at (3,5) {$F{\bm 2}$};
    \end{scope}
 \begin{scope} 
\draw[thick] (1)--(2);
\draw[thick] (1)--(3);
\draw[thick] (1)--(4);
\draw[thick] (1)--(5);
\draw[thick] (1)--(6);
\end{scope}
\end{tikzpicture}

(i) \(\mathcal{D}4\) requirements graph
\end{minipage}%
\begin{minipage}{.5\linewidth}
\centering
\begin{tikzpicture}[scale=0.8, transform shape]
\begin{scope} [vertex style/.style={draw,
                                       circle,
                                       minimum size=6mm,
                                       inner sep=0pt,
                                       outer sep=0pt}] 
\node[vertex style] (1) at (0,7) {$A{\bm 3}$};
\node[vertex style] (2) at (3,9) {$B{\bm 1}$};
\node[vertex style] (3) at (3,7) {$C{\bm 1}$};
\node[vertex style] (4) at (3,5) {$D{\bm 1}$};
\node[vertex style] (5) at (6,8) {$E{\bm 2}$};
\node[vertex style] (6) at (6,6) {$F{\bm 2}$};
    \end{scope}
 \begin{scope}
\draw[thick] (1)--(2);
\draw[thick] (1)--(3); 
\draw[thick] (1)--(4);
\draw[thick] (1)--(5);
\draw[thick] (1)--(6);
\draw[thick] (2)--(5);
\draw[thick] (2)--(6);
\draw[thick] (3)--(5);
\draw[thick] (3)--(6);
\draw[thick] (4)--(5);
\draw[thick] (4)--(6);
\end{scope}
\end{tikzpicture}  

(ii) \(\mathcal{D}4\) estimability graph
\end{minipage}%
\end{figure}
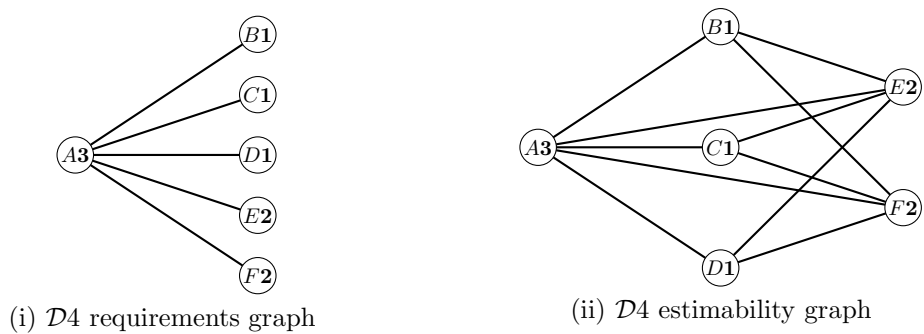

\subsection{Equivalent designs}
\label{iso}


Consider \(X_4\), the generator matrix for \(\mathcal{D}4\). This contains 
columns of \(\mathcal{X}_2\) with multiplicities one, two and three.  Any 
other  \(2^6\) factorial arranged in blocks of size four which also has a  generator matrix comprising columns of \(\mathcal{X}_2\) with multiplicities 
three, two and one will have an estimability graph in \(T(6,2)\) which is isomorphic to (or possibly the same as) the estimability graph for  
\(\mathcal{D}4\) and is  
said to be {\em equivalent} to \(\mathcal{D}4\). 
In general, two designs which have isomorphic estimability graphs in 
\(T(n,q)\) are equivalent in that they provide estimates of equivalent sets of two factor interactions. Note however that for \(q>2\) the designs themselves are not necessarily isomorphic: for example they may provide estimates of different numbers of three factor interactions. 

\vspace{0.1cm}

A design, \(\mathcal{D}\), comprising a \(2^n\) factorial in blocks of 
size \(2^q\) has generator matrix, \(X\), containing at least \(q\) 
independent members of \(\mathcal{X}_q\). Let the multiplicities of the 
members of \(\mathcal{X}_q\) in \(X\) be \(n_1 \ge n_2 \ge \cdots 
\ge n_{2^q-1}\ge 0\), where \(n_q>0\) and \(\sum_{i=1}^{2^q-1}n_i=n\). Then any other blocked \(2^n\) factorial 
having a full rank \(q \times n\) generator matrix with the same multiplicities is equivalent to 
\(\mathcal{D}\). 
Thus, subject to the requirement that the generator matrix has rank \(q\), the estimability properties of a blocked \(2^n\) factorial are wholly determined by the multiplicities of the  
generator matrix columns.  This prompts efficient representation of classes of designs with equivalent properties and of individual designs:


\vspace{0.2cm}

\noindent{\bf Definition:}~A \(2^n\) factorial in blocks of size \(2^q\) is said to have {\it profile set} \(\langle n_1, n_2, \ldots, n_{2^q-1} \rangle  \) 
where  \(n_1 \ge n_2 \ge \cdots \ge n_{2^q-1}\ge 0\) are the multiplicities 
of columns of \(\mathcal{X}_q\) in the generator matrix \(X\).

\vspace{0.2cm}

\noindent{\bf Definition:}~Consider a \(2^n\) factorial in blocks of size \(2^q\) with profile set \(\langle n_1, n_2, \ldots, n_{2^q-1} \rangle  \). Let \(q_0\) be the largest value such that \(n_{q_0}>0\). For \(i \le q_0\) denote the  factors 
relating to  entry \(n_{i}\)  in the profile set by \(F_{i,1}, \ldots, F_{i,n_{i}}\).
Then the design and its properties are determined by the {\it factor grouping} 
\(\{(F_{1,1}, \ldots, F_{1,n_{1}}),  \cdots, (F_{q_0,1}, \ldots, F_{q_0,n_{q_0}})\}\).

\vspace{0.1cm}

Two designs with the same profile set have isomorphic estimability graphs 
and are equivalent. The set of estimable interactions from a design are 
obtained directly from the factor grouping. This is illustrated by reference to 
\(\mathcal{D}3\) and \(\mathcal{D}4\).
Design \(\mathcal{D}4\) has profile   
set \(\langle 3,2,1 \rangle  \) and any  \(2^6\) factorial in blocks of size four with generator matrix comprising columns of \(\mathcal{X}_2\) with multiplicities one, two and three will have the same profile set and will be 
equivalent to \(\mathcal{D}4\). The estimability properties of \(\mathcal{D}
4\) are specified by the factor grouping \(\{(B,C,D), (E,F),(A)\}\): an 
interaction is estimable {\it iff} the factors are in different sets of the factor 
grouping: for example, \(BE\) is estimable but \(BC\) is not estimable. By comparison,  \(\mathcal{D}3\), with  profile  set  
\(\langle 2,2,2 \rangle  \) and factor grouping \(\{(A,B),(C,D),(E,F)\}\),  
is not equivalent to \(\mathcal{D}4\). Designs \(\mathcal{D}3\) and 
\(\mathcal{D}4\) are equivalent  to designs labelled 6-01/B4.1 and 
6-01/B4.2, considered in Case 1 of Sun {\it et al.}\,(1997), where they are 
specified in terms of four effects which generate the block defining 
contrast subgroup. (In fact, since \(q=2\), \(\mathcal{D}3\) is isomorphic to 
6-01/B4.1 and \(\mathcal{D}4\) is isomorphic to 6-01/B4.2). The design 
summaries in Sun {\it et al.}\,(1997) indicate that all main effects are 
estimable from both designs, and give the number of estimable interactions.  
In contrast with the generator matrix approach, the  
estimable interactions for each design are not evident without some work. 
Under H2 that effects of the same order are equally important, based on 
the number of estimable interactions 6-01/B4.1  (\(\mathcal{D}4\)) would 
be considered inferior to 6-01/B4.2  (\(\mathcal{D}3\)), which has minimum 
aberration blocking scheme. However, as demonstrated in Example 4, when 
H2 does not hold because a subset of interactions are of 
particular interest, there are situations in which 6-01/B4.1 
(\(\mathcal{D}4\)) would be preferred over 6-01/B4.2 (\(\mathcal{D}3\)).  

\subsection{Design construction via vertex colouring}\label{col}
In Figure \ref{eg4D4a}, {\it colours} \(\bm{ 1},\bm{2},\bm{3}\) are applied to vertices 
of the \(\mathcal{D}4\) requirements and estimability graphs, with vertices 
in the factor grouping set corresponding to \(n_i\) assigned colour \(\bm{i}\). For both graphs this results in a vertex colouring such that 
{\it adjacent} vertices, that is 
vertices which are joined by an edge, are coloured differently.  
In graph theory terms, the vertex colourings in Figure \ref{eg4D4a} constitute  proper {\it 3-colourings} and demonstrate that the 
graphs are {\it\(k\)-chromatic}, for some value \(k \le 3\). See Bondy and Murty (2008, Chapter 14).  If Example 4 was amended so that \(BD\), \(BE\) and \(DE\) were required as well as  interactions involving \(A\), this would give a 
new requirements graph with eight edges and, due to the arrangement of 
the edges,  it would not be possible to assign colours \(\bm{ 1},\bm{2},\bm{3}\) to the 
vertices to give a proper 3-colouring. Thus, there is no graph in \(T(6,2)\) with the new 
requirements graph as a subgraph and no \(2^6\) factorial in blocks of size four gives  estimates of all main effects and the eight interactions. The general point is summarised as follows:
\begin{lem}\label{lemcol}
The \(n\) main effects and the interactions of a requirements graph can be 
accommodated by a \(2^n\) factorial in blocks of size \(2^q\) {\it iff} the 
requirements graph  is \(k\)-chromatic for some value \(k \le 2^q-1\). 
\end{lem}
Lemma \ref{lemcol} prompts the strategy that will be adopted for design 
construction. A set of interactions is selected such that the requirements 
graph  is \(k\)-chromatic with \( k \le 2^q-1\). Let \(k_*=\min\{2^q-1,n\}\). A colouring in \(k_0\) 
colours is identified for some \(k_0\) such that \(\min\{q,k\}\le k_0 \le k_*\). The \(k_0\)-colouring corresponds to an estimability graph in \(T(n,q)\), with the requirements graph as a sub-graph, and to a 
factor grouping. The practitioner uses the \(k_0\)-colouring 
to identify a generator matrix for a \(2^n\) factorial in blocks of size \(2^q\) that yields all the required estimates. The factor grouping gives a concise means of recording the estimable interactions.

Note that if \(n\le 2^q-1\) then each vertex can be allocated a unique colour, to give an \(n\)-colouring. Any full rank matrix with columns comprising \(n\) members of \(\mathcal{X}_q\) will be a suitable generator matrix. The generator matrix \(X_2\) of Example 2 has this form and the corresponding estimability graph is a complete 6-partite graph in \(T(6,3)\). 
If \(n > 2^q-1\) then not all interactions can be estimated. Hence, there are sets of interactions which yield a requirements 
graph that is \(k\)-chromatic with \(k>2^q-1\). In such cases, no valid colouring is possible in \(2^q-1\) colours or fewer. Thus, the requirements graph will not be a subgraph of any graph in \(T(n,q)\) and no blocked \(2^n\) factorial will provide all the desired estimates. Known results in graph theory are useful in providing guidance on interaction sets that can be accommodated by designs. 

\begin{thm}\label{thmcondall}
Let \(\mathcal{S}\) denote a subset of  interactions, \(F_iF_j\), for an \(n\) factor process. If one or both of the following conditions are satisfied, then a  \(2^n\) factorial in blocks of size \(2^q\) exists from which all main effects and all interactions in 
\(\mathcal{S}\) can be estimated:
\begin{enumerate}
\item[(i)] \(\mathcal{S}\) contains no more than \(2^q-1\) interactions relating to any factor and there is no set of \(2^q\) factors for which \(\mathcal{S}\) contains all \(q\choose 2\) pairwise interactions; 
\item[(ii)] There are at most \(2^q-1\) factors for which \(\mathcal{S}\) contains at least \(2^q-1\) interactions.
\end{enumerate}
\end{thm}

\noindent{\bf Proof:}~Let \(\mathcal{G}\) be the requirements graph with 
\(n\) vertices and with edges relating to interactions in \(\mathcal{S}\). By Lemma \ref{lemcol}, a \(2^n\) design in blocks of size \(2^q\) can be obtained with the required properties {\it iff} \(\mathcal{G}\) is \(k\)-chromatic for some \(k \le 2^q-1\). The two cases are considered separately.
If \(\mathcal{S}\) satisfies (i), then the maximum vertex degree of \(\mathcal{G}\) does not exceed \(2^q-1\)  and \(\mathcal{G}\) does not contain the complete graph on \(2^q\) vertices.
By Brooks' Theorem (see Brooks (1941)), \(\mathcal{G}\) is \(k\)-chromatic for some \(k \le 2^q-1\). 
If \(\mathcal{S}\) satisfies (ii), then \(\mathcal{G}\) has fewer than \(2^q\) vertices of degree at least \(2^q-1\).
By Theorem 14.6 of Bondy and Murty (2008), every \(k\)-chromatic graph has at least \(k\) vertices of degree at least \(k-1\).  Thus \(\mathcal{G}\) is 
\(k\)-chromatic for some \(k \le 2^q-1\). The result follows.

\vspace{0.1cm}

When planning an experiment, if the practitioner lists interactions in order of priority then Theorem \ref{thmcondall} will give guidance on the selection of an interaction set which gives rise to a requirements graph with chromatic number no greater than \(2^q-1\). Obtaining a valid colouring with no more than \(2^q-1\) colours for such a graph is generally 
straightforward. However, if it proves challenging either to determine the chromatic number for a given requirements graph, or to obtain a valid colouring with no more than \(2^q-1\) colours, then graph theory routines in packages such as the Sage computer algebra system (see Stein {\it et al.}\,(2018))  are freely available.

\vspace{0.1cm}

\noindent{\bf Example~5  (part\,1):}~ {\it Two designs comprising \(2^7\) factorials in blocks of size four are sought, each with a different set of ten interactions of interest:
\begin{enumerate}
\item[] \(\mathcal{S}_1=\{AB,AC,AD,BC,BE,CD,DF,EF,EG,FG\}\);
\item[] \(\mathcal{S}_2=\{AB,AC,BC,BD, BE,CD,CF, CG,EF,EG\}\).
\end{enumerate}

Interaction sets \(\mathcal{S}_1\) and \(\mathcal{S}_2\) each satisfy exactly one of the conditions of Theorem \ref{thmcondall}, with 
\(\mathcal{S}_1\) satisfying (i) and \(\mathcal{S}_2\) satisfying (ii). The requirements graphs, complete with valid 
3-colourings, and the corresponding estimability graphs are presented in Figures \ref{eg31} and \ref{eg32}. 
Let \(\mathcal{DS}1\) and \(\mathcal{DS}2\)  denote designs that coincide with the graph colourings of Figures \ref{eg31} and \ref{eg32} and which accommodate \(\mathcal{S}_1\) and \(\mathcal{S}_2\). 
In both cases the factor grouping, and hence a generator matrix, can be written down from the requirements graph colouring: for example the factor grouping for \(\mathcal{DS}1\)  consistent with the 3-colouring in Figure \ref{eg31}\,(i) is \(\{(B,D,G),(A,E),(C,F)\}\) and a generator matrix for 
\(\mathcal{DS}1\)  is:
\begin{equation*}
X=\left( \begin{array}{ccccccc}
1 &1 &0& 1 &1& 0&1\\
0& 1&1&1&0 & 1 &1\end{array} \right).
\end{equation*}
Here, the  colours \({\bm 1},{\bm 2},{\bm 3}\) have been associated with columns of \(\mathcal{X}_2\) in the order \((1,1)^T,\) \((1,0)^T,\) \((0,1)^T\), 
but any one-to-one mapping between the colours (or factor sets) and the columns of \(\mathcal{X}_2\) is equivalent. Likewise, \(\mathcal{DS}2\) has factor grouping  \(\{(A,D,F),(B,G), (C,E)\}\).
Designs  \(\mathcal{DS}1\) and \(\mathcal{DS}2\) are equivalent: both 
have profile sets 
\(\langle 3,2,2 \rangle\) and isomorphic estimability graphs in \(T(7,2)\) and
give estimates of 16 interactions. }

\vspace{0.1cm}

Example 5 is returned to in subsequent sections.
\begin{figure}
\caption{Graphs for \(\mathcal{DS}1\) of Example 5 (part 1) with interaction set \(\mathcal{S}_1\)}\label{eg31}
\begin{minipage}{.4\linewidth}
    \centering
\vspace{0.2cm}
\hspace{0cm} \begin{tikzpicture}[scale=0.7, transform shape] 
   \begin{scope} [vertex style/.style={draw,
                                       circle,
                                       minimum size=6mm,
                                       inner sep=0pt,
                                       outer sep=0pt}] 
\node[vertex style] (1) at (0,6.5) {$A{\bm 2}$};
\node[vertex style] (2) at (2,8) {$B{\bm 1}$};
\node[vertex style] (3) at (2,6.5) {$C{\bm 3}$};
\node[vertex style] (4) at (2,5) {$D{\bm 1}$};
\node[vertex style] (5) at (4,8) {$E{\bm 2}$};
\node[vertex style] (6) at (4,5) {$F{\bm 3}$};
\node[vertex style] (7) at (6,6.5) {$G{\bm 1}$};
    \end{scope}
 \begin{scope} 
\draw[thick] (1)--(2);
\draw[thick] (1)--(3);
\draw[thick] (1)--(4);
\draw[thick] (2)--(3);
\draw[thick] (3)--(4);
\draw[thick] (2)--(5);
\draw[thick] (4)--(6);
\draw[thick] (5)--(6);
\draw[thick] (5)--(7);
\draw[thick] (6)--(7);
\end{scope}
\end{tikzpicture}

Requirements graph
\end{minipage}%
~~
\begin{minipage}{.4\linewidth}
\centering
\vspace{0.2cm}
\begin{tikzpicture}[scale=0.7, transform shape] 
   \begin{scope} [vertex style/.style={draw,
                                       circle,
                                       minimum size=6mm,
                                       inner sep=0pt,
                                       outer sep=0pt}] 
\node[vertex style] (2) at (4.5,9) {$D{\bm 1}$};
\node[vertex style] (5) at (7.8,8) {$E{\bm 2}$};
\node[vertex style] (7) at (12.5,7.7) {$F{\bm 3}$};
\node[vertex style] (6) at (12.5,10.3) {$C{\bm 3}$};
\node[vertex style] (4) at (7.8,10) {$A{\bm 2}$};
\node[vertex style] (3) at (4.5,7.4) {$G{\bm 1}$};
\node[vertex style] (1) at (4.5,10.6) {$B{\bm 1}$};
    \end{scope}
 \begin{scope} 
\draw[thick] (1)--(4);
\draw[thick] (1)--(5);
\draw[thick] (1)--(6);
\draw[thick] (1)--(7);
\draw[thick] (2)--(4);
\draw[thick] (2)--(5);
\draw[thick] (2)--(6);
\draw[thick] (2)--(7);
\draw[thick] (3)--(4);
\draw[thick] (3)--(5);
\draw[thick] (3)--(6);
\draw[thick] (3)--(7);
\draw[thick] (4)--(6);
\draw[thick] (4)--(7);
\draw[thick] (5)--(6);
\draw[thick] (5)--(7);
\end{scope}
\end{tikzpicture}

Estimability graph
\end{minipage}%
\end{figure}
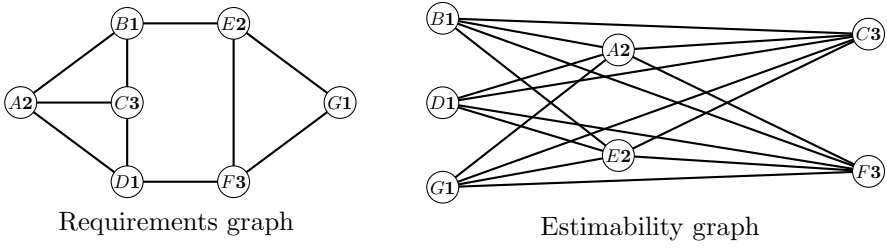
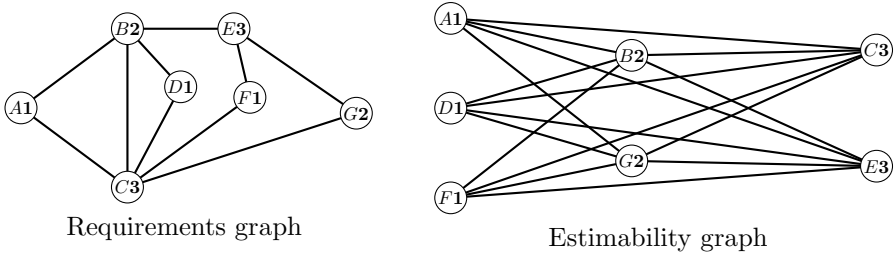
\begin{figure}
\caption{Graphs for \(\mathcal{DS}2\) of Example 5 (part 1) with interaction set \(\mathcal{S}_2\)}\label{eg32}
\begin{minipage}{.4\linewidth}
    \centering
\vspace{0.2cm}
\hspace{0cm} \begin{tikzpicture}[scale=0.7, transform shape] 
   \begin{scope} [vertex style/.style={draw,
                                       circle,
                                       minimum size=6mm,
                                       inner sep=0pt,
                                       outer sep=0pt}] 
\node[vertex style] (1) at (0,6.5) {$A{\bm 1}$};
\node[vertex style] (2) at (2,8) {$B{\bm 2}$};
\node[vertex style] (3) at (2,5) {$C{\bm 3}$};
\node[vertex style] (4) at (3,6.9) {$D{\bm 1}$};
\node[vertex style] (5) at (4,8) {$E{\bm 3}$};
\node[vertex style] (6) at (4.3,6.7) {$F{\bm 1}$};
\node[vertex style] (7) at (6.3,6.4) {$G{\bm 2}$};
    \end{scope}
 \begin{scope} 
\draw[thick] (1)--(2);
\draw[thick] (1)--(3);
\draw[thick] (2)--(3);
\draw[thick] (3)--(4);
\draw[thick] (2)--(4);
\draw[thick] (2)--(5);
\draw[thick] (3)--(6);
\draw[thick] (5)--(7);
\draw[thick] (5)--(6);
\draw[thick] (3)--(7);
\end{scope}
\end{tikzpicture}

Requirements graph
\end{minipage}%
~~
\begin{minipage}{.4\linewidth}
\centering
\vspace{0.2cm}
\begin{tikzpicture}[scale=0.7, transform shape] 
   \begin{scope} [vertex style/.style={draw,
                                       circle,
                                       minimum size=6mm,
                                       inner sep=0pt,
                                       outer sep=0pt}] 
\node[vertex style] (2) at (4.5,9) {$D{\bm 1}$};
\node[vertex style] (5) at (7.9,8) {$G{\bm 2}$};
\node[vertex style] (7) at (12.5,7.9) {$E{\bm 3}$};
\node[vertex style] (6) at (12.5,10.1) {$C{\bm 3}$};
\node[vertex style] (4) at (7.9,10) {$B{\bm 2}$};
\node[vertex style] (3) at (4.5,7.3) {$F{\bm 1}$};
\node[vertex style] (1) at (4.5,10.7) {$A{\bm 1}$};
    \end{scope}
 \begin{scope} 
\draw[thick] (1)--(4);
\draw[thick] (1)--(5);
\draw[thick] (1)--(6);
\draw[thick] (1)--(7);
\draw[thick] (2)--(4);
\draw[thick] (2)--(5);
\draw[thick] (2)--(6);
\draw[thick] (2)--(7);
\draw[thick] (3)--(4);
\draw[thick] (3)--(5);
\draw[thick] (3)--(6);
\draw[thick] (3)--(7);
\draw[thick] (4)--(6);
\draw[thick] (4)--(7);
\draw[thick] (5)--(6);
\draw[thick] (5)--(7);
\end{scope}
\end{tikzpicture}

Estimability graph
\end{minipage}%
\end{figure}


\subsection{Alternative colourings and profile sets}\label{eq}

A graph which is \(k\)-chromatic may admit  multiple proper \(k_0\)-colourings, for some \(k_0 \ge k\), and these will not necessarily share the same 
set cardinalities. 
Thus, a requirements graph with \(k \le k_*\) may permit different \(k_0\)-colourings, for some \(k_0\) with \(\max\{k,q\} \le k_0 \le k_*\), leading to blocked \(2^n\) factorials with different profile sets 
and therefore able to estimate different numbers of interactions. 
Further, 
for a requirements graph that is \(k\)-chromatic with \(k <k_*\), a colouring in \(k_*\) colours can always be found.  These points are now investigated briefly in the context of constructing designs with desirable properties, the 
assumption being that if two alternative \(2^n\) factorials in blocks of size \(2^q\) with different profile sets both give the required estimates, then a design with profile set resulting in a larger number of estimable interactions is preferred. The first point is now demonstrated.

\vspace{0.1cm}

\noindent{\bf Example 5 (part 2):}~{\it The 3-colouring of the requirements 
graph for \(\mathcal{S}_2\) given in Figure \ref{eg32} yields \(\mathcal{DS}2\) with profile set \(\langle 3,2,2 \rangle\). In Figure 
\ref{eg52}, the \(\mathcal{S}_2\) requirements graph is reproduced with  
an alternative 3-colouring. This gives design \(\mathcal{DS}2A\) with profile set \(\langle 3,3,1 \rangle\). The estimability graph for \(\mathcal{DS}2A\) is 
included in 
Figure \ref{eg52}. Of \(\mathcal{DS}2\) and \(\mathcal{DS}2A\), the former is preferred since it gives estimates of 16 interactions compared to the 15 interactions estimable from \(\mathcal{DS}2A\). 
A 3-colouring of the \(\mathcal{S}_2\)  requirements graph can also be given leading to a design with profile set  \(\langle 4,2,1 \rangle\), yielding estimates of only 14 interactions. }
\begin{figure}
\caption{Graphs for \(\mathcal{DS}2A\) of Example 5 (part 2) with interaction set \(\mathcal{S}_2\)}\label{eg52}
\begin{minipage}{.4\linewidth}
    \centering
\vspace{0.2cm}
\hspace{0cm} \begin{tikzpicture}[scale=0.7, transform shape] 
   \begin{scope} [vertex style/.style={draw,
                                       circle,
                                       minimum size=6mm,
                                       inner sep=0pt,
                                       outer sep=0pt}] 
\node[vertex style] (1) at (0,6.5) {$A{\bm 1}$};
\node[vertex style] (2) at (2,8) {$B{\bm 2}$};
\node[vertex style] (3) at (2,5) {$C{\bm 3}$};
\node[vertex style] (4) at (3,6.9) {$D{\bm 1}$};
\node[vertex style] (5) at (4,8) {$E{\bm 1}$};
\node[vertex style] (6) at (4.3,6.7) {$F{\bm 2}$};
\node[vertex style] (7) at (6.3,6.4) {$G{\bm 2}$};
    \end{scope}
 \begin{scope} 
\draw[thick] (1)--(2);
\draw[thick] (1)--(3);
\draw[thick] (2)--(3);
\draw[thick] (3)--(4);
\draw[thick] (2)--(4);
\draw[thick] (2)--(5);
\draw[thick] (3)--(6);
\draw[thick] (5)--(7);
\draw[thick] (5)--(6);
\draw[thick] (3)--(7);
\end{scope}
\end{tikzpicture}

Requirements graph
\end{minipage}%
~~
\begin{minipage}{.4\linewidth}
\centering
\vspace{0.2cm}
\begin{tikzpicture}[scale=0.7, transform shape] 
   \begin{scope} [vertex style/.style={draw,
                                       circle,
                                       minimum size=6mm,
                                       inner sep=0pt,
                                       outer sep=0pt}] 
\node[vertex style] (2) at (4.5,9.7) {$D{\bm 1}$};
\node[vertex style] (5) at (12.5,9.7) {$F{\bm 2}$};
\node[vertex style] (7) at (4.5,8.7) {$E{\bm 1}$};
\node[vertex style] (6) at (8.5,7) {$C{\bm 3}$};
\node[vertex style] (4) at (12.5,10.7) {$B{\bm 2}$};
\node[vertex style] (3) at (12.5,8.7) {$G{\bm 2}$};
\node[vertex style] (1) at (4.5,10.7) {$A{\bm 1}$};
    \end{scope}
 \begin{scope} 
\draw[thick] (1)--(3);
\draw[thick] (1)--(4);
\draw[thick] (1)--(5);
\draw[thick] (1)--(6);
\draw[thick] (2)--(3);
\draw[thick] (2)--(4);
\draw[thick] (2)--(5);
\draw[thick] (2)--(6);
\draw[thick] (3)--(6);
\draw[thick] (3)--(7);
\draw[thick] (4)--(6);
\draw[thick] (4)--(7);
\draw[thick] (5)--(6);
\draw[thick] (5)--(7);
\draw[thick] (6)--(7);
\end{scope}
\end{tikzpicture}

Estimability graph
\end{minipage}%
\end{figure}

\vspace{0.1cm}

To investigate the second point, consider a \(k\)-chromatic  requirements graph, where \(k<k_*\), together with a \(k_0\)-colouring for some \(k_0\) such that \(k_0 <k_*\). 
If \(k_0 \ge q\) this colouring corresponds to a design with profile set  
having 
\(n_1>1\), \(n_i\ge 1\), for \(i=2, \ldots, k_0\), and \(n_i=0\), for \(i=k_0+1, \ldots,2^q-1\). 
A blocked \(2^n\) factorial with full rank generator matrix comprising 
copies of \(k_0\) columns of \(\mathcal{X}_q\), arranged according to the profile set provides the required estimates. 
The graph colouring can be adjusted to 
incorporate additional colours. 
To demonstrate this, arbitrarily select \([n_1/2] \)  vertices of colour \({\bm 1}\)
and recolour these with colour \(k_0+1\). The new colouring is a valid 
\((k_0+1)\)-colouring
 which  leads to an estimability graph in \(T(n,q)\) and a blocked \(2^n\) factorial 
giving estimates of an additional \([n_1/2](n_1-[n_1/2])\) 
interactions over those achieved by the original design. 
The process can be repeated until a  colouring is achieved with \(k_*\) colours.
If \(k_0 <q\), the initial \(k_0\)-colouring will not correspond to a blocked \(2^n\) 
factorial, since the colouring does not yield a full rank generator matrix. Use of the process as described for \(k_0 \ge q\) will result in  graph colourings using progressively more colours. Once a \(q\)-colouring is achieved, the graph corresponds to a blocked \(2^n\) factorial  and, after successively including more colours, ultimately a  \(k_*\)-colouring is achieved.

In terms of the profile set, the number of estimable interactions from a blocked \(2^n\) factorial is \(\sum_{i=1}^{2^q-2}\sum_{j=i+1}^{2^q-1}n_in_j\), which achieves the bound \(\phi_{\max}\) of Theorem  \ref{int} 
when \(n_1-n_{2^q-1} \le 1\). In the notation of Theorem 
\ref{int} this occurs {\it iff} both \( n_1 \le v+1\) and \(n_{2^q-1}=v\).
Such a colouring is termed 
an {\it equitable} \(k_*\)-colouring.
If \(n \le 2^q-1\) an equitable \(n\)-colouring has each vertex coloured uniquely: design \(\mathcal{D}2\) of Example 2, with profile set \(\langle 1,1,1,1,1,1,0 \rangle\), is such a case. 
The colourings in Figure \ref{eg31} and \ref{eg32}, leading to \(\mathcal{DS}1\) and \(\mathcal{DS}2\) are equitable \(k_*\)-colourings, with \(k_*=3\). 
For \(n > 2^q-1\), not all 
\(k\)-chromatic graphs with \(k \le 2^q-1\) admit equitable \(k_*\)-colourings: the requirements graph of \(\mathcal{D}3\), displayed in Figure \ref{eg4D4}(i) cannot be equitably 
3-coloured. 

\section{Blocked \(2^{n-p}\) factorial designs}

Unless the number of factors is small it is unlikely that resources will be available to conduct a full blocked \(2^n\) factorial. Instead, a \(2^{n-p}\) 
factorial is chosen, where \(p>0\), and the treatment combinations of the principal fraction are arranged in blocks of size \(2^q\) according to a generating 
matrix to form a blocked \(2^{n-p}\) factorial. A \(2^{n-p}\)  factorial is specified by \(p\) independent defining words
\(F_1^{z_{i,1}}\cdots F_{n-p}^{z_{i,n-p}}F_i\) for \(i=n-p+1, \ldots, n\),  with \(z_{i,j} \in \{0,1\}\), summarised by the \(p \times (n-p) \) matrix 
\(Z=(z_{i,j})\). The defining words  
give rise to a treatment defining contrast sub-group containing \(2^{p}-1\) effects. The factorial effects not in the treatment defining contrast sub-group are partitioned into \(2^{n-p}-1\) alias sets, each of size \(2^p\).

The complexity of planning an experiment   to estimate all main effects and selected  interactions is far greater when using a blocked 
\(2^{n-p}\)  factorial than when using a blocked 
\(2^{n}\)  factorial. This is because consideration needs to be given to both the principal block generating 
sub-group and the treatment defining contrast sub-group. The major difference between the constructions of blocked \(2^n\) factorials and blocked \(2^{n-p}\) factorials proposed here is that, for 
the former the practitioner is able to select all \(n\) columns of \(X\) whereas, for the latter, only the first \(n-p\) columns are selected. 
A blocked \(2^{n-p}\) factorial has generator matrix \(X=(X^I X^{II})\). The columns of the \(q \times (n-p)\) sub-matrix \(X^{I}\) are selected from \(\mathcal{X}_q\) and the \(q \times p\) sub-matrix \(X^{II}\) is then determined from \(X^{II}=X^I Z^T\), where calculations are modulo 2. 
In a \(2^{n-p}\) factorial arranged in blocks of size \(2^q\), of the alias sets, \(2^{n-p-q}-1\) are confounded with blocks and the remaining \(2^{n-p}-2^{n-p-q}\) are orthogonal to blocks. Whether a factorial effect is confounded with blocks or not is determined directly from \(X\) as in \S\ref{full}. 
Theorem \ref{thm1} is adjusted to include fractional designs:
\begin{thm}\label{ffthm1}
For a \(2^{n-p}\) factorial in blocks of size \(2^q\), an effect is not confounded with blocks {\it iff} at least one of the generating treatment combinations contains an odd number of effect terms.
\end{thm}
Likewise, for effects of interest:
\begin{cor}\label{ffcor}
For a \(2^{n-p}\) factorial in blocks of size \(2^q\), all main effects are unconfounded with blocks {\it iff} every column of \(X\)  is a column of \(\mathcal{X}_q\).
\end{cor}
\begin{cor}\label{ffcor2}
For a \(2^{n-p}\) factorial in blocks of size \(2^q\),  an interaction is not confounded with blocks {\it iff} the corresponding two columns of \(X\) are different.
\end{cor}
An effect of interest will be estimable {\it iff} both conditions are satisfied:
\begin{enumerate}
\item[C1:] all \(2^p-1\) aliases of the effect can be considered negligible and can be suppressed;
\item[C2:] the effect is not confounded with blocks.
\end{enumerate} 
Since all main effects are required and no two factor interactions are assumed negligible, only \(2^{n-p}\) factorials of resolution at least IV 
are used. Similarly, only selections of columns from \(\mathcal{X}_q\) for \(X^I\) which yield \(X^{II}\) sub-matrices also with columns all from \(\mathcal{X}_q\), that is with no zero sum columns, are entertained. 
The design construction approach is illustrated:

\vspace{0.1cm}

\noindent{\bf Example 6:}~{\it A design comprising a \(2^{6-1}\) factorial 
arranged in blocks of size eight is required, to provide estimates of all main effects and as many interactions as possible. 

A sensible starting point is to use a fraction with defining word of length five 
or six to yield a \(2^{6-1}\) factorial for which
main effects and two factor interactions are only aliased with negligible 
effects. Consider the defining word \(ABCDEF\), giving 
\(Z=(1,1,1,1,1)\). We seek a \(3 \times 6\) generator matrix 
\(X=(X^I X^{II})\) by selecting five of the seven members of \(\mathcal{X}_3\) for the 
columns of \(X^I\) and then obtaining \(X^{II}\) using \(X^{II}=X^I Z^T\). However, any choice of five distinct members of \(\mathcal{X}_3\)  for 
\(X^I\) yields a vector \(X^{II}\) which is one of the selected members. Thus the generator matrix contains a repeated column. For example:
\begin{equation*}
X=\left( \begin{array}{cccccc}
1&0&0&1&1&1 \\
0&1&1&1&1&0\\
0&0&1&0&1&0 \end{array} \right).
\end{equation*}
With this generator matrix, \(AF\) is confounded with blocks and so is not estimable, but all main effects and the other 14 interactions are estimable. Clearly the design can be adjusted so that the interaction thought to be of least interest is non-estimable.

The same duplicated column property is observed if a defining word of length five, such as \(ABCDF\), is used. A defining word of length four and careful choice of members of \(\mathcal{X}_3\) for \(X^{I}\) can yield a matrix \(X\) with six distinct columns. However, although the resulting design will have no main effects or two factor interactions confounded with blocks, three pairs of two factor interactions will be aliased and so fail C2.
  }
  
\vspace{0.1cm}

\noindent{\bf Example 7:}~{\it An experiment involving an industrial process is planned in which a 
\(2^{9-2}\) factorial, in factors \(A,B, \ldots, I\),  is arranged in blocks of size eight. Due to their roles in the process, it is thought that \(G,H\) and \(I\) do not interact with each other, but interactions between these factors and the other six factors are of interest. Thus a design with profile set \(\langle 3,1,1,1,1,1,1 \rangle\), with \(G,H\) and \(I\) forming the factor group of size three, would provide the required estimates and also estimates of interactions between \(A, \ldots, F\).
Use of the defining words \(ABEGH\) and \(ABCDEFI\) gives a resolution V fraction and has 
\begin{equation*}
Z=\left( \begin{array}{ccccccc}
1&1&0&0&1&0&1 \\
1&1&1&1&1&1&0\end{array} \right).
\end{equation*} 
Then, trial and error in selecting and arranging six different columns of \(\mathcal{X}_3\) in \(X^I\) yields the generator matrix:
\begin{equation*}
X=\left( \begin{array}{ccccccccc}
1&0&0&0&1&1&1&1&1 \\
0&1&0&1&1&0&1&1&1\\
0&0&1&1&0&1&1&1&1 \end{array} \right),
\end{equation*}
which provides the estimates summarised above.}

\section{Blocks of size four}
\label{four}

In this section we focus on the construction of \(2^n\) and \(2^{n-p}\) 
designs in blocks of size four. For a factorial experiment arranged in blocks 
of size two, a design which provides estimates of all main effects will not 
enable estimation of any  interactions. Thus, the factorial experiments  in blocks 
of size four is the class with smallest block size enabling estimation of all main effects and some interactions. Further,
the blocks are small enough for blocking to be effective in most cases.
With the large number of factorial effects confounded with blocks, designs 
in blocks of size four with desirable properties are challenging to construct. The importance of these designs and the challenge of their  construction merit the special attention given to the design class.
To arrange a \(2^n\) factorial in blocks of size four  we use the approach of 
\S\ref{col} to obtain information on properties of interaction sets which give rise to requirements graphs with chromatic number 2 or 3. In the following result, conditions (i) and (ii) are covered by Theorem \ref{thmcondall} but are included here for completeness.
\begin{thm}\label{thmcond}
Let \(\mathcal{S}\) denote a subset of  interactions, \(F_iF_j\), for an \(n\) factor process. If one or more of the following are satisfied, then a  \(2^n\) factorial in blocks of size four exists from which all main effects and all interactions in 
\(\mathcal{S}\) can be estimated:
\begin{enumerate}
\item[(i)] \(\mathcal{S}\) contains no more than three interactions relating to any factor and there is no set of four factors for which \(\mathcal{S}\) contains all six pairwise interactions;
\item[(ii)] There are at most three factors for which \(\mathcal{S}\) contains more than two interactions;
\item[(iii)] There is no subset of factors, \(F_i\), for \(i=1, \ldots, m\ge 3\), for which all \(m\) interactions \(F_1F_2, F_2F_3, \ldots F_{m-1}F_m,F_mF_1 \in \mathcal{S}\), or there is at least one factor that is common to all such subsets. 
\end{enumerate}
\end{thm}

\noindent{\bf Proof of (iii):}~Let \(\mathcal{G}\) be the requirements graph with \(n\) vertices and with edges relating to interactions in \(\mathcal{S}\). By Lemma \ref{lemcol}, a \(2^n\) design in blocks of size four can be obtained with the required properties {\it iff} \(\mathcal{G}\) is \(k\)-chromatic for \(k \le 3\). 

If \(\mathcal{S}\) satisfies (iii) by there being no subset of factors, \(F_1,\ldots,F_m\), for \(m\ge 3\), such that \(F_1F_2, F_2F_3, \ldots F_{m-1}F_m,F_mF_1 \in \mathcal{S}\) then \(\mathcal{G}\) 
contains no cycles, that is \(\mathcal{G}\)  is {\it acyclic}, and is therefore 2-chromatic, (see Bondy and Murty (2008, Chapters 4,\,14)). Suppose instead that 
\(\mathcal{S}\) 
satisfies (iii) by there being at least one subset of factors, \(F_1,\ldots,F_m\), for \(m\ge 3\), forming a cycle in \(\mathcal{G}\), with all such subsets having at least one common factor. Without loss of generality let 
\(F_{1}\) be a common factor for all subsets. Removal of \(F_{1}\) and adjacent edges from \(\mathcal{G}\) gives a graph \(\mathcal{G}'\) in \(n-1\) vertices which is acyclic 
and is therefore 2-chromatic. Consider any 2-colouring of the vertices of \(\mathcal{G}'\) and reinstate \(F_{1}\) and incident edges to reform 
\(\mathcal{G}\). Vertex \(F_{1}\) can be coloured with a third colour to give a valid 3-colouring of \(\mathcal{G}\). Hence \(\mathcal{G}\) is \(k\)-chromatic, with \(k \le 3\).

\vspace{0.1cm}

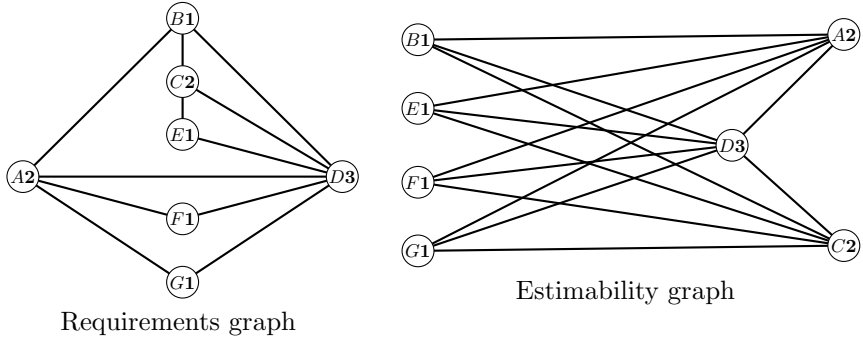
\begin{figure}
\caption{Graphs for \(\mathcal{DS}3\) of Example 5 (part 3) with interaction set \(\mathcal{S}_3\)}\label{eg33}
\begin{minipage}{.4\linewidth}
    \centering
\vspace{0.2cm}
\hspace{0cm} \begin{tikzpicture}[scale=0.7, transform shape] 
   \begin{scope} [vertex style/.style={draw,
                                       circle,
                                       minimum size=6mm,
                                       inner sep=0pt,
                                       outer sep=0pt}] 
\node[vertex style] (1) at (-1,5) {$A{\bm 2}$};
\node[vertex style] (2) at (2,8) {$B{\bm 1}$};
\node[vertex style] (3) at (2,6.8) {$C{\bm 2}$};
\node[vertex style] (4) at (2,5.8) {$E{\bm 1}$};
\node[vertex style] (5) at (2,4.2) {$F{\bm 1}$};
\node[vertex style] (6) at (2,3) {$G{\bm 1}$};
\node[vertex style] (7) at (5,5) {$D{\bm 3}$};
    \end{scope}
 \begin{scope} 
\draw[thick] (1)--(2);
\draw[thick] (1)--(5);
\draw[thick] (1)--(6);
\draw[thick] (1)--(7);
\draw[thick] (4)--(3);
\draw[thick] (3)--(7);
\draw[thick] (4)--(7);
\draw[thick] (5)--(7);
\draw[thick] (7)--(6);
\draw[thick] (2)--(3);
\draw[thick] (2)--(7);
\end{scope}
\end{tikzpicture}

Requirements graph
\end{minipage}%
\begin{minipage}{.4\linewidth}
\centering
\begin{tikzpicture}[scale=0.7, transform shape] 
   \begin{scope} [vertex style/.style={draw,
                                       circle,
                                       minimum size=6mm,
                                       inner sep=0pt,
                                       outer sep=0pt}] 
\node[vertex style] (2) at (4.5,8.8) {$F{\bm 1}$};
\node[vertex style] (5) at (12.5,11.6) {$A{\bm 2}$};
\node[vertex style] (7) at (12.5,7.6) {$C{\bm 2}$};
\node[vertex style] (6) at (10.4,9.5) {$D{\bm 3}$};
\node[vertex style] (4) at (4.5,10.2) {$E{\bm 1}$};
\node[vertex style] (3) at (4.5,7.5) {$G{\bm 1}$};
\node[vertex style] (1) at (4.5,11.5) {$B{\bm 1}$};
    \end{scope}
 \begin{scope} 
\draw[thick] (1)--(5);
\draw[thick] (1)--(6);
\draw[thick] (1)--(7);
\draw[thick] (2)--(5);
\draw[thick] (2)--(6);
\draw[thick] (2)--(7);
\draw[thick] (3)--(5);
\draw[thick] (3)--(6);
\draw[thick] (3)--(7);
\draw[thick] (4)--(5);
\draw[thick] (4)--(6);
\draw[thick] (4)--(7);
\draw[thick] (5)--(6);
\draw[thick] (6)--(7);
\end{scope}
\end{tikzpicture}

Estimability graph
\end{minipage}%
\end{figure}

Use of Theorem \ref{thmcond} is demonstrated by the following continuation of Example 5: 

\vspace{0.1cm}

\noindent{\bf Example~5 (part\,3):}~{\it Designs comprising \(2^7\) replicates in blocks of size four are sought, with the following interaction sets of interest:
\begin{enumerate}
\item[]  \(\mathcal{S}_3=\{AB,AD,AF,AG,BC,BD,CD,CE,DE,DF,DG\}\);
\item[]  \(\mathcal{S}_4=\{AB,AC,AD,AE,AG, BF, CD,CG,DG, EF\}\).
\end{enumerate}

Set  \(\mathcal{S}_3\) satisfies condition (iii) of Theorem \ref{thmcond}. 
The requirements graph, complete with valid 
3-colouring, and the corresponding estimability graph are presented in Figure \ref{eg33}. Let \(\mathcal{DS}3\)  denote the design that coincides with the graph colouring and thus accommodates \(\mathcal{S}_3\).
With profile set  \(\langle 4,2,1 \rangle\),  \(\mathcal{DS}3\) gives estimates of 14 interactions. The factor grouping from the colouring is \(\{(B,E,F,G),(A,C),(D)\}\).
Set \(\mathcal{S}_4\) does not satisfy any of the conditions of Theorem \ref{thmcond}. The requirements graph for 
\(\mathcal{S}_4\), displayed in 
Figure \ref{eg34},  is 4-chromatic. No \(2^7\) factorial in blocks of size four provides estimates of all main effects and all interactions in 
\(\mathcal{S}_4\). Theorem \ref{thmcond} can be used to give guidance on subsets of \(\mathcal{S}_4\) that can be accommodated.
The removal of any interaction involving at least one of \(C,D,G\) gives an interaction set that satisfies one or both  of conditions (ii) and (iii). }

\vspace{0.1cm}

Sets \(\mathcal{S}_1-\mathcal{S}_3\) each satisfy exactly one condition of 
Theorem \ref{thmcond}. 
It is interesting to consider estimability properties of the two   \(2^7\) factorials in blocks of size four,  labelled 7-0.1/B5.1 and 7-0.1/B5.2 in Sun {\it et al.}\,(1997). After some work it can be shown that these designs have profile sets \(\langle 3,2,2 \rangle\) and \(\langle 3,3,1 \rangle\) and provide  
estimates of all main effects and of 16 and 15 interactions, respectively. Thus \(\mathcal{DS}1\) and \(\mathcal{DS}2\), of Example 4 (part 1), are isomorphic to  7-0.1/B5.1. Design,  \(\mathcal{DS}3\) is not isomorphic to either 7-0.1/B5.1 or 7-0.1/B5.2. Further, neither of these designs can accommodate interaction set \(\mathcal{S}_3\), even with a relabelling of factors, since the requirements graph for \(\mathcal{S}_3\) does not admit a 3-colouring consistent with profile sets \(\langle 3,2,2 \rangle\) or 
\(\langle 3,3,1 \rangle\).
\begin{center}
\begin{figure}
\caption{Requirements graph for interaction set \(\mathcal{S}_4\) of Example 5 (part 3)}\label{eg34}
\vspace{0.2cm}
\hspace{2cm} \begin{tikzpicture}[scale=0.7, transform shape] 
   \begin{scope} [vertex style/.style={draw,
                                       circle,
                                       minimum size=6mm,
                                       inner sep=0pt,
                                       outer sep=0pt}] 
\node[vertex style] (1) at (0,6.5) {$A$};
\node[vertex style] (2) at (1.6,8.2) {$B$};
\node[vertex style] (3) at (1.6,4.8) {$C$};
\node[vertex style] (4) at (4,8.5) {$D$};
\node[vertex style] (5) at (4,4.5) {$E$};
\node[vertex style] (6) at (6,7.5) {$F$};
\node[vertex style] (7) at (6,5.5) {$G$};
    \end{scope}
 \begin{scope} 
\draw[thick] (1)--(2);
\draw[thick] (1)--(3);
\draw[thick] (1)--(4);
\draw[thick] (1)--(5);
\draw[thick] (1)--(7);
\draw[thick] (2)--(6);
\draw[thick] (3)--(4);
\draw[thick] (3)--(7);
\draw[thick] (4)--(7);
\draw[thick] (5)--(6);
\end{scope}
\end{tikzpicture}
\end{figure}
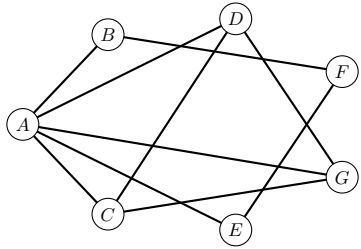
\end{center}

\vspace{-0.7cm}

Every acyclic graph with at least one edge is 2-chromatic and therefore 3-colourable. A result of Bollob\'as and Guy (1983) establishes that many acyclic 
graphs are equitably 3-colourable. The result is reproduced without proof, expressed in terms of an interaction set for an \(n\) factor process, rather than in 
the graph theory terms of Bollob\'as and Guy: 

\begin{thm}\label{thmtree}{\bf Bollob\'as and Guy}
Let \(\mathcal{S}\) denote a subset of interactions for an \(n\) factor process such that there is no subset of factors \(F_i\), for \(i=1, \ldots, m\ge 3\), for which all \(m\) interactions \(F_1F_2, F_2F_3, \ldots F_{m-1}F_m,F_mF_1 \in \mathcal{S}\). Then a design comprising a \(2^n\) factorial in blocks of size four which provides the required estimates and achieves the bound \(\phi_{\max}\) of Theorem \ref{int} exists if one of the following conditions is satisfied:
\begin{enumerate}
\item[(i)] \(\mathcal{S}\) contains no more than \((n+8)/3\) interactions involving any one factor;
\item[(ii)] The largest number of interactions involving any factor is \((n+10)/3\).
\end{enumerate} 
\end{thm} 
\begin{cor}\label{cortree}
Consider an \(n\) factor process and interaction set  \(\mathcal{S}\) with property as specified in Theorem \ref{thmtree}. If \(\mathcal{S}\) contains no more 
than four interactions involving any factor then a  \(2^n\) factorial in blocks of size four exists which yields the required estimates and achieves the bound of Theorem \ref{int}.  
\end{cor}
\noindent{\bf Proof:}~Let \(n_0\) be the largest number of interactions in \(\mathcal{S}\) involving a common factor. Then \(n_0 \le n-1\). For \(n_0 \le 4\), 
we have \(n_0 < (n_0+9)/3 \le (n+8)/3\) and condition\,(i) of Theorem \ref{thmtree} holds.

\vspace{0.1cm}

For a situation with an acyclic requirements graph, Theorem \ref{thmtree} and Corollary \ref{cortree} give simple conditions guaranteeing the existence of a design which provides estimates of all effects of interest and achieves the bound \(\phi_{\max}\). 

\vspace{0.1cm}

\noindent{\bf Example~8:}~{\it A \(2^6\) replicate in blocks of size four is sought to provide estimates of all interactions in  \(\mathcal{S}=\{AB,AC,AD,AE, EF\}\).

\vspace{0.1cm}

Set \(\mathcal{S}\) satisfies Corollary \ref{cortree} and thus the  requirements graph (not included)  has an equitable 3-colouring. An example of such a 
colouring leads to the  factor grouping \(\{(A,F),(B,C),(D,E)\}\). Thus, a 
generator matrix that yields a \(2^6\) factorial in blocks of size four which 
provides the required estimates and achieves \(\phi_{\max}\) is:
\begin{equation*}
X=\left( \begin{array}{cccccc}
1 &1 &1& 0&0& 1\\
0& 1&1&1&1 & 0\end{array} \right).
\end{equation*}}
\subsection{Designs in blocks of size four from \(2^{n-p}\) factorials of resolution at least V}\label{secres5}
We consider the construction of blocked \(2^{n-p}\) factorials based on \(2^{n-p}\) factorials of resolution at least V. This has the advantage that  
both main effects and interactions will only be aliased with negligible effects. Therefore, C1 is satisfied for all interactions and therefore any interaction of interest that is not confounded with blocks will be estimable. For given \((n,p)\) pair for which \(2^{n-p}\) factorials 
of resolution at least V exist, the approach is to work systematically through \(2^{n-p}\) factorials of appropriate resolution in order of aberration, starting with a fractional factorial of minimum aberration, to identify blocked \(2^{n-p}\) factorials covering as many profile sets as possible. In this way we 
build up a portfolio of design templates. Existing catalogues are a source of \(2^{n-p}\) factorials in order of aberration and details of \(2^{n-p}\) factorials
used are contained in Table \ref{tab:fg}  in the Appendix. The process is demonstrated by the following example.

\vspace{0.1cm}

\noindent{\bf Example~9:}~{\it Design templates are sought for \(2^{7-1}\) factorials in blocks of size four to yield estimates of all main effects and a subset of interactions. 

\vspace{0.1cm}

For \(n=7\), the profile sets with \(n_3>0\) are  \(\langle 3,2,2 \rangle\), \(\langle 3,3,1 \rangle\), \(\langle 4,2,1 \rangle\) and \(\langle 5,1,1 \rangle\). The \(2^{7-1}\) 
 factorial with minimum aberration is 7-1.1 with fraction generator \(F_1F_2F_3F_4F_5F_6F_7\). Thus \(Z=(1,1,1,1,1,1)\). Potential generator matrices are constructed as 
follows: without loss of generality \((1,0)^T\) is assigned to the first column of \(X^I\). All possible \(3^5\) allocations of columns of \(\mathcal{X}_2\) to columns \(2-6\) 
of \(X^I\) are considered. For each allocation, the \(2 \times 1\) sub-matrix \(X^{II}\) is obtained as the sum of the columns of \(X^I\), working modulo 2. An 
allocation that yields a column of \(\mathcal{X}_2\) as \(X^{II}\) signifies a blocked \(2^{7-1}\) replicate from which all main effects are estimable. The factor grouping is stored if the profile set has 
not already been identified. The scanning exercise is undertaken using Matlab R2017a. The fractional factorial  7-1.1 yields generator 
matrices with profile sets \(\langle 5,1,1 \rangle\) and \(\langle 3,3,1 \rangle\). The factor groupings recorded are 
\(\{ (F_1,F_2,F_3,F_4,F_5),(F_6),(F_7)\}\) and \(\{(F_2,F_3,F_4),(F_5, F_6,F_7),(F_1)\}\).
To accommodate other profile sets, attention moves to the fractional factorial 7-1.2, which has fraction generator \(F_1F_2F_3F_4F_5F_7\), giving  \(Z=(1,1,1,1,0,1)\). This yields generator 
matrices with  \(X^{II}\) having non-zero sum, for profile sets \(\langle 4,2,1 \rangle\) and \(\langle 3,2,2 \rangle\). The factor groupings recorded are \(\{(F_2,F_3,F_4,F_5),(F_1,F_7),\) \( (F_6)\}\) and \(\{(F_2, F_3,F_6),(F_1, F_7),(F_4,F_5)\}\).
No further \(2^{7-1}\)  factorials are considered since a generator matrix has now been found for every profile set with \(n_3>0\). Any profile set with only two positive terms can be subsumed into one of these sets. }

\vspace{0.1cm}

It is informative to return again to Example 5 to illustrate use of the design templates.

\vspace{0.1cm}

\noindent{\bf Example~5\,(part\,4):}~{\it Blocked \(2^{7-1}\) factorials, with blocks of size four, are required to provide estimates of interactions in sets \(\mathcal{S}_1\), \(\mathcal{S}_2\) and \(\mathcal{S}_3\) of Example 5 (parts 1,3).

\vspace{0.1cm}

Denote the planned designs by \(\mathcal{DS}11\), \(\mathcal{DS}21\) and \(\mathcal{DS}31\). The graphs in Figures \ref{eg31},\,\ref{eg32} and \ref{eg33} from 
Example 5 (parts 1,3) still apply. As with \(\mathcal{DS}1\), design \(\mathcal{DS}11\) can be constructed with profile set \(\langle 3,2,2 \rangle\) and factor grouping 
\(\{(B,D,G),(A,E),(C,F)\}\). A mapping of \(A, \ldots,G\) onto \(F_1, \ldots,F_7\) which is consistent with the factor grouping recorded for profile set 
\(\langle 3,2,2 \rangle\) in Example 8 is: 
\(A \rightarrow F_1, \,E\rightarrow F_7,\,C \rightarrow F_4,\,F \rightarrow F_5,\,B \rightarrow F_2,\,D \rightarrow F_3,\, G \rightarrow F_6\), and the  
defining word is \(ABCDEF\). Design \(\mathcal{DS}21\) is also constructed with profile set \(\langle 3,2,2 \rangle\). With factor grouping \(\{(A,D,F),(B,G),(C,E)\}\), using Example 9 a suitable mapping  
  is: \(B \rightarrow F_1,\,G \rightarrow F_7,\,C \rightarrow F_4,\, E \rightarrow F_5,\,A \rightarrow F_2, \,D\rightarrow F_3,\,F \rightarrow F_6\) and the defining word is \(ABCDEG\).  Design \(\mathcal{DS}31\) has 
profile set   \(\langle 4,2,1 \rangle\) and factor grouping \(\{(B,E,F,G),(A,C),(D)\}\). A mapping  consistent with Example 9 is: \(D \rightarrow F_6,\,A \rightarrow F_1, \,C\rightarrow F_7,\,B \rightarrow F_2,\,E \rightarrow F_3,\,F \rightarrow F_4,\,G \rightarrow F_5\) and the defining word is 
\(ABCEFG\).

The designs are constructed from knowledge of the factor grouping and defining word. For example, the treatment combinations in 
\(\mathcal{DS}31\)  are those in the principal \(2^{7-1}\) factorial with defining word \(ABCEFG\). A suitable generator matrix is:
\begin{equation*}
X=
\left( \begin{array}{ccccccc}
1 &0 & 1 &1& 0&0 &0\\
0& 1&0&1 & 1 &1&1\end{array} \right),
\end{equation*}
giving principal block \((1), acd, bdefg, abcefg\). The remaining 15 blocks contain the other treatment combinations of the half replicate, with each block containing a coset of the treatment combinations of the principal block.}

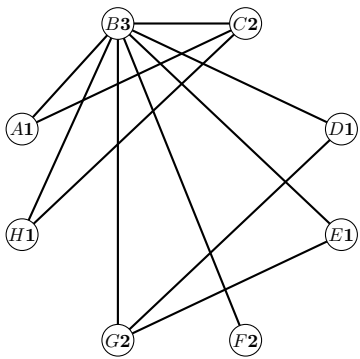
\begin{figure}
\caption{Requirements graph for Example 10}\label{eg8}
    \centering
\vspace{0.2cm}
\hspace{2cm} \begin{tikzpicture}[scale=0.7, transform shape] 
   \begin{scope} [vertex style/.style={draw,
                                       circle,
                                       minimum size=6mm,
                                       inner sep=0pt,
                                       outer sep=0pt}] 
\node[vertex style] (1) at (0,6) {$A{\bm 1}$};
\node[vertex style] (2) at (1.8,8) {$B{\bm 3}$};
\node[vertex style] (3) at (4.2,8) {$C{\bm 2}$};
\node[vertex style] (4) at (6,6) {$D{\bm 1}$};
\node[vertex style] (5) at (6,4) {$E{\bm 1}$};
\node[vertex style] (6) at (4.2,2) {$F{\bm 2}$};
\node[vertex style] (7) at (1.8,2) {$G{\bm 2}$};
\node[vertex style] (8) at (0,4) {$H{\bm 1}$};
    \end{scope}
 \begin{scope} 
\draw[thick] (2)--(1);
\draw[thick] (2)--(3);
\draw[thick] (2)--(4);
\draw[thick] (2)--(5);
\draw[thick] (2)--(6);
\draw[thick] (2)--(7);
\draw[thick] (2)--(8);
\draw[thick] (1)--(3);
\draw[thick] (3)--(8);
\draw[thick] (4)--(7);
\draw[thick] (5)--(7);
\end{scope}
\end{tikzpicture}
\end{figure}

For \(5 \le n \le 11\), Tables \ref{tab:res5} and \ref{tab:res5a} give design templates for blocked \(2^{n-p}\) factorials with up to 128 runs using 
\(2^{n-p}\) factorials with resolution at least V. 
Table \ref{tab:res5} encompasses all profile sets for \((n,p)\) pairs \((6,1),(7,1)\) and \((8,2)\), and Table \ref{tab:res5a} covers all 
profile sets for pairs \((8,1)\) and \((9,2)\).      
For pair \((5,1)\),  Table \ref{tab:res5} gives  a design template for the  profile set \(\langle 3,1,1 \rangle\). No \(2^{5-1}\) factorial in blocks of size four 
exists with the same estimability capability as Example 1, which has profile set \(\langle 2,2,1 \rangle\).  For pair \((10,3)\), Table \ref{tab:res5a}  
gives templates for all profile sets except \(\langle 8,1,1 \rangle\), and for pair \((11,4)\), templates 
are given for all profile sets with every factor group containing at least two members. Use of the tables is demonstrated:

\vspace{0.1cm}

\noindent{\bf Example~10:}~{\it  A \(2^{8-2}\) factorial arranged in blocks of size four is required to assess all interactions involving factor \(B\) and interactions \(AC,CH,DG, EG\).

\vspace{0.1cm}

Since \(2^{8-2}\) factorials in blocks of size four are given for all profile sets in Table \ref{tab:res5}, the problem is similar to that of arranging a full 
\(2^8\) factorial in blocks of size four.  
The set of interactions satisfies conditions (ii) and (iii) of Theorem \ref{thmcond} and so can be represented by a requirements graph which is 3-chromatic. 
We start by producing the requirements graph (see Figure \ref{eg8}) and investigating possible 3-colourings. Since all interactions involving \(B\) are 
required, the colour used for vertex \(B\) cannot be used for any other vertex.
Thus, \(n_3=1\). Of the profile sets having this property, \(\langle 4,3,1 \rangle\) gives a design 
with the largest number of estimable interactions. There are several colourings which result in profile set \(\langle 4,3,1 \rangle\). One such colouring  is 
included in Figure \ref{eg8} and has factor grouping  \(\{(A,D,E,H),(C,F,G),(B)\}\). 
Using Table \ref{tab:res5}, a blocked \(2^{8-2}\) factorial can be constructed via mapping 
\(B \rightarrow F_1,\,C \rightarrow F_2, \,F\rightarrow F_6,\,G \rightarrow F_8,\,A \rightarrow F_3,\,D \rightarrow F_4,\,E \rightarrow F_5,\,H \rightarrow F_7\) and defining words \(ABCDH\) and \(BCEFG\). A generator matrix is 
\begin{equation*}
X=
\left( \begin{array}{cccccccc}
1 &0 & 1 &1& 1&1 &1&1\\
0& 1&1&0 & 0 &1&1&0\end{array} \right),
\end{equation*}
giving principal block \((1), acdefgh, bcfg, abdeh\). The remaining 15 blocks are constructed as cosets to give a design in 64 runs. }

\begin{table}
\centering
\begin{minipage}{11cm}
\caption{Templates for \(2^{n-p}\) factorials in blocks of size four with up to 64 runs, from \(2^{n-p}\) factorials of resolution at least V}\label{tab:res5}

\begin{tabular}{ccccccc}\hline
\(n\) & \(p\)& runs & profile  & int\(^a\) & fractional & factor \\
 &  & &set &   & factorial & grouping \\ \hline
5 & 1 &16& \(\langle 3,1,1 \rangle\) & 7 &   5-1.1 & \(\{(F_1,F_2,F_3),(F_4),(F_5)\}\) \\
6 & 1 & 32 &\(\langle 2,2,2 \rangle\) & 12 &   6-1.1 & \(\{ (F_1,F_6),(F_2,F_3),(F_4,F_5)\}\) \\
6 & 1 & 32 &\(\langle 3,2,1 \rangle\) & 11 &   6-1.2 & \(\{ (F_2,F_3,F_6),(F_4,F_5),(F_1)\}\) \\
6 & 1 & 32 &\(\langle 4,1,1 \rangle\) & 9 &   6-1.2 & \(\{ (F_1,F_2,F_3,F_5),(F_4),(F_6)\}\) \\
7 & 1 & 64&\(\langle 3,2,2 \rangle\) & 16 &   7-1.2 & \(\{ (F_2,F_3,F_6),(F_1,F_7),(F_4,F_5)\}\) \\
7 & 1 &64& \(\langle 3,3,1 \rangle\) & 15 &   7-1.1 & \(\{ (F_2,F_3,F_4),(F_5,F_6,F_7),(F_1)\}\) \\
7 & 1 &64& \(\langle 4,2,1 \rangle\) & 14 &   7-1.2 & \(\{ (F_2,F_3,F_4,F_5),(F_1,F_7),(F_6)\}\) \\
7 & 1 &64& \(\langle 5,1,1 \rangle\) & 11 &   7-1.1 & \(\{ (F_1,F_2,F_3,F_4,F_5),(F_6),(F_7)\}\) \\
8 & 2 &64& \(\langle 3,3,2 \rangle\) & 21 &   8-2.1 & \(\{ (F_1,F_4,F_7),(F_2,F_5,F_6),(F_3,F_8)\}\) \\
8 & 2 &64& \(\langle 4,2,2 \rangle\) & 20 &   8-2.1 & \(\{ (F_1,F_2,F_3,F_5),(F_4,F_6),(F_7,F_8)\}\) \\
8 & 2 &64& \(\langle 4,3,1 \rangle\) & 19 &   8-2.1 & \(\{ (F_3,F_4,F_5,F_7), (F_2,F_6,F_8),(F_1)\}\) \\
8 & 2 &64& \(\langle 5,2,1 \rangle\) & 17 &   8-2.1 & \(\{ (F_2,F_4,F_5,F_6,F_7),(F_3,F_8),(F_1)\}\) \\
8 & 2 &64& \(\langle 6,1,1 \rangle\) & 13 &   8-2.1 & \(\{ (F_3,F_4,F_5,F_6,F_7,F_8),(F_1),(F_2)\}\) \\
\hline
\end{tabular}

\vspace{0.1cm}

\(^a\)number of estimable interactions
\end{minipage}
\end{table}

\begin{table}
\centering
\begin{minipage}{12cm}
\caption{Templates for \(2^{n-p}\) factorials in blocks of size four with 128 runs, from \(2^{n-p}\) factorials of resolution at least V~~} \label{tab:res5a}

\begin{tabular}{cccccc}\hline
\(n\) & \(p\)&  profile  & int\(^a\)  & fractional & factor  \\ 
 &  & set &   & factorial & grouping \\ \hline
8 & 1 & \(\langle 3,3,2 \rangle\) & 21 &   8-1.2 & \(\{ (F_2,F_3,F_4),(F_5,F_6,F_8), (F_1,F_7)\}\) \\
8 & 1 & \(\langle 4,2,2 \rangle\) & 20 &   8-1.1 & \(\{ (F_2,F_3,F_4,F_5),(F_1,F_8),(F_6,F_7)\}\) \\
8 & 1 & \(\langle 4,3,1 \rangle\) & 19 &   8-1.2 & \(\{(F_2,F_3,F_4,F_7),(F_5,F_6,F_8), (F_1)\}\) \\
8 & 1 & \(\langle 5,2,1 \rangle\) & 17 &   8-1.2 & \(\{ (F_2,F_3,F_4,F_5,F_8),(F_6,F_7),(F_1)\}\) \\
8 & 1 & \(\langle 6,1,1 \rangle\) & 13 &   8-1.2 & \(\{ (F_1,F_2,F_3,F_4,F_5,F_7),(F_6),(F_8)\}\) \\
9 & 2 & \(\langle 3,3,3 \rangle\) & 27 &   9-2.1 & \(\{ (F_1,F_8,F_9),(F_2,F_4,F_6),(F_3,F_5,F_7)\}\) \\
9 & 2 & \(\langle 4,3,2 \rangle\) & 26 &   9-2.1 & \(\{ (F_2,F_3,F_6,F_9),(F_1,F_7,F_9),(F_4,F_5)\}\) \\
9 & 2 & \(\langle 5,2,2 \rangle\) & 24 &   9-2.1 & \(\{ (F_1,F_6,F_7,F_8,F_9),(F_2,F_3),(F_4,F_5)\}\) \\
9 & 2 & \(\langle 4,4,1 \rangle\) & 24 &   9-2.2 & \(\{ (F_1,F_4,F_8,F_9),(F_2,F_5,F_6,F_7),(F_3)\}\) \\
9 & 2 & \(\langle 5,3,1 \rangle\) & 23 &   9-2.2 & \(\{ (F_3,F_5,F_6,F_7,F_9),(F_1,F_2,F_8),(F_4)\}\) \\
9 & 2 & \(\langle 6,2,1 \rangle\) & 20 &   9-2.2 & \(\{ (F_2,F_4,F_5,F_6,F_7,F_8),(F_1,F_9),(F_3)\}\) \\
9 & 2 & \(\langle 7,1,1 \rangle\) & 15 &   9-2.2 & \(\{ (F_1,F_2,F_3,F_5,F_6,F_7,F_9),(F_4),(F_8)\}\) \\
10 & 3 & \(\langle 4,3,3 \rangle\) & 33 &  10-3.1 & \(\{ (F_3,F_5,F_7,F_{10}),(F_1,F_8,F_9),(F_2,F_4,F_6)\}\) \\
10 & 3 & \(\langle 4,4,2 \rangle\) & 32 &  10-3.1 & \(\{ (F_1,F_6,F_8,F_{10}),(F_2,F_3,F_7,F_9),(F_4,F_5)\}\) \\
10 & 3 & \(\langle 5,3,2 \rangle\) & 31 &  10-3.1 & \(\{ (F_2,F_3,F_6,F_9,F_{10}),(F_1,F_7,F_8),(F_4,F_5)\}\) \\
10 & 3 & \(\langle 5,4,1 \rangle\) & 29 &  10-3.1 & \(\{ (F_2,F_4,F_6,F_7,F_9),(F_1,F_3,F_5,F_8),(F_{10})\}\) \\
10 & 3 & \(\langle 6,2,2 \rangle\) & 28 &  10-3.1 & \(\{ (F_1,F_6,F_7,F_8,F_9,F_{10}),(F_2,F_3),(F_4,F_5)\}\)\\
10 & 3 & \(\langle 6,3,1 \rangle\) & 27 &  10-3.1 & \(\{ (F_2,F_3,F_4,F_6,F_7,F_9),(F_1,F_5,F_9),(F_{10})\}\) \\
10 & 3 & \(\langle 7,2,1 \rangle\) & 23 &  10-3.1 & \(\{ (F_2,F_4,F_5,F_6,F_7,F_8,F_9),(F_1,F_3),(F_{10})\}\) \\
11 & 4 & \(\langle 4,4,3 \rangle\) & 40 &  11-4.1 & \(\{ (F_2,F_4,F_6,F_{11}),(F_3,F_5,F_7,F_{10}),(F_1,F_8,F_9)\}\) \\
11 & 4 & \(\langle 5,3,3 \rangle\) & 39 &  11-4.1 & \(\{ (F_2,F_4,F_6,F_7,F_9),(F_1,F_3,F_{11}),(F_5,F_8,F_{10})\}\) \\
11 & 4 & \(\langle 5,4,2 \rangle\) & 38 &  11-4.1 & \(\{ (F_2,F_3,F_6,F_9,F_{10}),(F_1,F_7,F_8,F_{11}),(F_4,F_5)\}\) \\
11 & 4 & \(\langle 6,3,2 \rangle\) & 36 &  11-4.1 & \(\{ (F_4,F_5,F_6,F_7,F_{10}, F_{11}),(F_1,F_8,F_9),(F_2,F_3)\}\) \\
11 & 4 & \(\langle 7,2,2 \rangle\) & 32 &  11-4.1 & \(\{ (F_1,F_6,F_7,F_8,F_9,F_{10}, F_{11}),(F_2,F_3),(F_4,F_5)\}\) \\
\hline
\end{tabular}

\vspace{0.1cm}
   
\(^a\)number of estimable interactions
\end{minipage}
\end{table}

\subsection{Designs in blocks of size four from \(2^{n-p}\) factorials of resolution IV}\label{secres4}

We give brief consideration to the construction of blocked \(2^{n-p}\) factorials based on fractional factorials of resolution IV. As in \S\ref{secres5}, main 
effects will only be aliased with effects which are negligible. However, some two factor interactions will be aliased with each other. 
To be specific, each four factor interaction in the  treatment defining contrast sub-group will correspond to three pairs of aliased two factor interactions. We 
focus on blocked \(2^{n-p}\) factorials obtained from \(2^{n-p}\) factorials of resolution IV for which the treatment defining contrast subgroup contains only one word of length four. 
Since the columns of \(X\) pertaining to this four factor interaction necessarily sum to \((0,0)^T\), the four columns are either four copies of the same 
column of \(\mathcal{X}_2\) or two columns of \(\mathcal{X}_2\) each copied twice. Thus, the four constituent factors of the interaction will all be contained in one 
factor group or will be distributed between two factor groups with each containing a pair. In the former situation, the factors of the aliased two factor interactions occur in the same factor group. Thus, the interactions fail C2 and so also failing C1 has no impact on the estimability capability of the design. In the latter case, four interactions involving factors in two factor groups will be inestimable: these will fail only C1.  

Table \ref{tab:res4a} gives design templates for \(n=7,9,12\) not covered in Tables \ref{tab:res5} and \ref{tab:res5a}. Where a design can be found for an 
\((n,p)\) pair and profile set combination with all aliased two-factor interactions involving factors in the same factor group then such a design is the only one 
recorded, 
since it improves on designs with aliased interactions spread between two factor groups. See for example the blocked \(2^{7-2}\) factorial with profile set 
\(\langle 1,1,5 \rangle\) in the third row of Table \ref{tab:res4a}. 
Designs with common \((n,p)\) pair and profile set combination with aliased interactions spread between two factor groups do not necessarily have equivalent 
estimability structure.  For example, consider the two blocked \(2^{9-3}\) factorials both with profile set \(\langle 4,3,2 \rangle\) in the seventh and eighth 
rows of Table \ref{tab:res4a}. Both have the factors of aliased interactions split between two factor groups but the estimability graphs are not isomorphic. 
Both estimability graphs are formed from the same complete 3-partite graph, with four edges deleted. 
In the first case the deleted edges are between the vertex sets (factor groups) of sizes two and four, whereas in the second they are between the vertex 
sets (factor groups) of sizes three and four.

Designs in which factors involved in aliased two-factor interactions are all contained in the same factor group are used in exactly the same way as those in Tables \ref{tab:res5} and 
\ref{tab:res5a}. For an \((n,p)\) pair and profile set  with aliased two-factor interactions involving two factor groups, care has to be given to 
allocating factors within factor groups to take account of the non-estimable interactions, where possible. 
To illustrate the consequence of alternative graph colourings, and of different mappings within factor groups for a given colouring, we return one final time to Example 5.

\vspace{0.1cm}

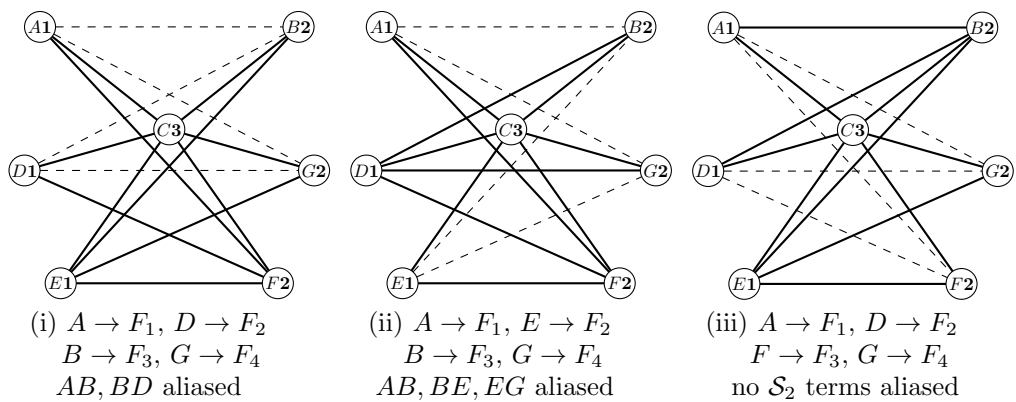
\begin{figure}
\caption{Estimability graphs for blocked \(2^{7-2}\) factorials with interaction set \(\mathcal{S}_2\)}\label{eg32again}
\begin{minipage}{.25\linewidth}
 \centering
\vspace{0.2cm}
\begin{tikzpicture}[scale=0.68, transform shape] 
   \begin{scope} [vertex style/.style={draw,
                                       circle,
                                       minimum size=6mm,
                                       inner sep=0pt,
                                       outer sep=0pt}] 
\node[vertex style] (1) at (0,8) {$A{\bm 1}$};
\node[vertex style] (2) at (5,8) {$B{\bm 2}$};
\node[vertex style] (4) at (-0.3,5.2) {$D{\bm 1}$};
\node[vertex style] (7) at (5.3,5.2) {$G{\bm 2}$};
\node[vertex style] (6) at (4.6,3) {$F{\bm 2}$};
\node[vertex style] (5) at (0.4,3) {$E{\bm 1}$};
\node[vertex style] (3) at (2.5,6) {$C{\bm 3}$};
    \end{scope}
 \begin{scope} 
\draw[thick] (1)--(3);
\draw[thick] (2)--(3);
\draw[thick] (3)--(4);
\draw[thick] (3)--(5);
\draw[thick] (3)--(6);
\draw[thick] (3)--(7);
\draw[dashed,-] (1)--(2);
\draw[thick] (1)--(6);
\draw[dashed,-] (1)--(7);
\draw[dashed,-] (4)--(2);
\draw[thick] (4)--(6);
\draw[dashed,-] (4)--(7);
\draw[thick] (5)--(2);
\draw[thick] (5)--(6);
\draw[thick] (5)--(7);
\end{scope}
\end{tikzpicture}

(i) \(A \rightarrow F_1,\,D \rightarrow F_2\)

~~~\(B \rightarrow F_3,\, G \rightarrow F_4\)

\(AB,BD\) aliased
\end{minipage}%
~~~~~~
\begin{minipage}{.25\linewidth}
 \centering
\vspace{0.2cm}
\begin{tikzpicture}[scale=0.68, transform shape] 
   \begin{scope} [vertex style/.style={draw,
                                       circle,
                                       minimum size=6mm,
                                       inner sep=0pt,
                                       outer sep=0pt}] 
\node[vertex style] (1) at (0,8) {$A{\bm 1}$};
\node[vertex style] (2) at (5,8) {$B{\bm 2}$};
\node[vertex style] (4) at (-0.3,5.2) {$D{\bm 1}$};
\node[vertex style] (7) at (5.3,5.2) {$G{\bm 2}$};
\node[vertex style] (6) at (4.6,3) {$F{\bm 2}$};
\node[vertex style] (5) at (0.4,3) {$E{\bm 1}$};
\node[vertex style] (3) at (2.5,6) {$C{\bm 3}$};
    \end{scope}
 \begin{scope} 
\draw[thick] (1)--(3);
\draw[thick] (2)--(3);
\draw[thick] (3)--(4);
\draw[thick] (3)--(5);
\draw[thick] (3)--(6);
\draw[thick] (3)--(7);
\draw[dashed,-] (1)--(2);
\draw[thick] (1)--(6);
\draw[dashed,-] (1)--(7);
\draw[thick] (4)--(2);
\draw[thick] (4)--(6);
\draw[thick] (4)--(7);
\draw[dashed,-] (5)--(2);
\draw[thick] (5)--(6);
\draw[dashed,-] (5)--(7);
\end{scope}
\end{tikzpicture}

(ii) \(A \rightarrow F_1,\,E \rightarrow F_2\)

~~~\(B \rightarrow F_3,\, G \rightarrow F_4\)

~\(AB,BE,EG\) aliased
\end{minipage}%
~~~~~~
\begin{minipage}{.25\linewidth}
 \centering
\vspace{0.2cm}
\begin{tikzpicture}[scale=0.68, transform shape] 
   \begin{scope} [vertex style/.style={draw,
                                       circle,
                                       minimum size=6mm,
                                       inner sep=0pt,
                                       outer sep=0pt}] 
\node[vertex style] (1) at (0,8) {$A{\bm 1}$};
\node[vertex style] (2) at (5,8) {$B{\bm 2}$};
\node[vertex style] (4) at (-0.3,5.2) {$D{\bm 1}$};
\node[vertex style] (7) at (5.3,5.2) {$G{\bm 2}$};
\node[vertex style] (6) at (4.6,3) {$F{\bm 2}$};
\node[vertex style] (5) at (0.4,3) {$E{\bm 1}$};
\node[vertex style] (3) at (2.5,6) {$C{\bm 3}$};
    \end{scope}
 \begin{scope} 
\draw[thick] (1)--(3);
\draw[thick] (2)--(3);
\draw[thick] (3)--(4);
\draw[thick] (3)--(5);
\draw[thick] (3)--(6);
\draw[thick] (3)--(7);
\draw[thick] (1)--(2);
\draw[dashed,-] (1)--(6);
\draw[dashed,-] (1)--(7);
\draw[thick] (4)--(2);
\draw[dashed,-] (4)--(6);
\draw[dashed,-] (4)--(7);
\draw[thick] (5)--(2);
\draw[thick] (5)--(6);
\draw[thick] (5)--(7);
\end{scope}
\end{tikzpicture}

(iii) \(A \rightarrow F_1,\,D \rightarrow F_2\)

~~~~\(F \rightarrow F_3,\, G \rightarrow F_4\)

~~~no \(\mathcal{S}_2\) terms aliased
\end{minipage}%
\end{figure}

\noindent{\bf Example 5 (part\,5):}~{\it A design comprising a \(2^{7-2}\) factorial in blocks of size four is sought to provide estimates of interactions in set \(\mathcal{S}_2\).

\vspace{0.1cm}

The requirements graph and colouring in Figure \ref{eg32} correspond to profile set \(\langle 3,2,2 \rangle\) and factor grouping 
\(\{(A,D,F), (B,G),(C,E)\}\). The top row of Table  \ref{tab:res4a} suggests a mapping: 
\(B \rightarrow F_1,\,G \rightarrow F_6,\,C \rightarrow F_2,\, E \rightarrow F_3,\,A \rightarrow F_4, \,D\rightarrow F_5,\,F \rightarrow F_7\) with 
defining words \(BCEG\) and \(ABCDF\). However, four of the selected interactions, namely \(BC,BE,CG\) and \(EG\), are aliased with other two-factor interactions and so are not 
estimable. Bar a reordering of the colour labels, Figure \ref{eg32} gives the only colouring with 
profile \(\langle 3,2,2 \rangle\). Thus, the 32 run design will provide estimates of all main effects but of only six of the interactions in \(\mathcal{S}_2\). We 
consider colourings leading to alternative profile sets. Figure \ref{eg52} gives a colouring for the requirements graph of \(\mathcal{S}_2\) consistent with 
profile set  \(\langle 3,3,1 \rangle\), and factor grouping \(\{(A,D,E),(B,F,G), (C)\}\). The construction from the second row of Table \ref{tab:res4a} will give a suitable design but with interactions involving factors mapped onto \(F_1-F_4\) aliased. To construct such a design, two factors from the factor group comprising \(A,D,E\) are mapped onto  \(F_1,F_2\) and 
two factors from \(B,F,G\) are mapped onto onto \(F_3,F_4\). The estimability graphs for three mappings are illustrated in Figure 
\ref{eg32again}.  Each consists of a complete 3-partite graph consistent with profile set \(\langle 3,3,1 \rangle\) with four edges removed. The edges to be removed relate to interactions that are 
inestimable despite involving factors from different factor groups and are indicated by dotted lines for clarity. The mappings in 
Figure \ref{eg32again}\,(i) and (ii) lead to two and three  interactions of \(\mathcal{S}_2\) being inestimable. The final mapping of Figure \ref{eg32again}\,(iii) is preferred over these since it 
provides all the required estimates. The full mapping of  Figure \ref{eg32again}\,(iii)  is  \(C \rightarrow F_5,\,A \rightarrow F_1,\,D \rightarrow F_2,\,E \rightarrow F_7\), \(B \rightarrow F_6,\,F \rightarrow F_3,\, G \rightarrow F_4\). A \(2^{7-2}\) factorial in blocks of size four consistent with this has generator matrix:
\begin{equation*}
X=
\left( \begin{array}{ccccccc}
1 &0 & 1 &1& 1&0 &0\\
0& 1&1&0 & 0 &1&1\end{array} \right),
\end{equation*}
and defining words \(ABEF\) and \(ACDEG\).
This shows that a profile set able to give estimates of fewer  interactions can be advantageous and that different mappings result in different numbers of interactions of \(\mathcal{S}_2\) being estimable.}

\begin{table}
\centering
\begin{minipage}{16cm}
\caption{Templates for blocked \(2^{n-p}\) factorials with up to 128 runs, from \(2^{n-p}\) factorials of resolution IV} \label{tab:res4a}
\resizebox{\textwidth}{!}{%
\begin{tabular}{cccccclc}\hline
\(n\) & \(p\)&runs&  profile  & int\(^a\)   & fractional & ~~~~~~~~~~~~~~~~~factor & aliased\(^b\)   \\
 &  && set &   & factorial & ~~~~~~~~~~~~~~~grouping & interactions\\ \hline
 7 & 2 & 32& \(\langle 3,2,2 \rangle\) & 12 &   7-2.1 &\(\{(F_4,F_5,F_7),(F_1,F_6),(F_2,F_3)\}\)&\(F_1F_2,F_1F_3,F_2F_6,F_3F_6\)\\
7 & 2 &32&\(\langle 3,3,1 \rangle\) & 11 &   7-2.1 & \(\{(F_1,F_2,F_7), (F_3,F_4,F_6),(F_5)\}\) &\(F_1F_3,F_1F_6,F_2F_3,F_2F_6\)\\
7 & 2 & 32&\(\langle 5,1,1 \rangle\) & 11 &   7-2.1 & \(\{(F_1,F_2,F_3,F_4,,F_7),(F_5),(F_7)\}\) & none\\
7 & 2 & 32&\(\langle 4,2,1 \rangle\) & 10 &   7-2.1 & \(\{(F_1,F_5,F_6,F_7),(F_2,F_3),(F_4)\}\) &\(F_1F_2,F_1F_3,F_2F_6,F_3F_6\)\\
9 & 3 & 64&\(\langle 5,2,2 \rangle\) & 24 &   9-3.1 & \(\{(F_1,F_2,F_3,F_4,F_7), (F_5,F_6),(F_8,F_9) \}\) &none\\
9 & 3 & 64&\(\langle 3,3,3 \rangle\) & 23 &   9-3.1 & \(\{(F_1,F_3,F_9),(F_2,F_6,F_7), (F_4,F_5,F_8)\}\) &\(F_1F_2,F_1F_7,F_2F_3,F_3F_7\)\\
9 & 3 & 64&\(\langle 4,3,2 \rangle\) & 22 &   9-3.1 & \(\{ (F_1,F_5,F_7,F_8),(F_4,F_6,F_9),(F_2,F_3)\}\) &\(F_1F_2,F_1F_3,F_2F_7,F_3F_7\)\\
9 & 3 & 64&\(\langle 4,3,2 \rangle\) & 22 &   9-3.1 & \(\{ (F_2,F_4,F_5,F_7),(F_1,F_3,F_9),(F_6,F_8)\}\) &\(F_1F_2,F_1F_7,F_2F_3,F_3F_7\)\\
9 & 3 & 64&\(\langle 6,2,1 \rangle\) & 20 &   9-3.1 & \(\{(F_1,F_2,F_3,F_7,F_8,F_9),(F_5,F_6),(F_4)\}\) &none\\
12 & 5 &128& \(\langle 4,4,4 \rangle\) & 48 &   12-5.1 & \(\{(F_1,F_8,F_9,F_{12}),(F_2,F_4,F_6,F_{11}),\) &none\\
&&&&&&~~~~~~~~\((F_3,F_5,F_7,F_{10})\}\)&\\
12 & 5 & 128&\(\langle 5,4,3 \rangle\) & 47 &   12-5.2 & \(\{(F_3,F_4,F_6,F_{10},F_{12}),(F_2,F_5,F_7,F_{11}),\) &none\\
&&&&&&~~~~~~~~\( (F_1,F_8,F_9)\}\)&\\
12 & 5 & 128&\(\langle 6,4,2 \rangle\) & 44 &   12-5.1 & \(\{(F_4,F_5,F_6,F_7,F_{10},F_{11}), (F_1,F_8,F_9,F_{12}),\) &none\\
&&&&&&~~~~~~~~\( (F_2,F_3)  \}\)&\\
12 & 5 &128& \(\langle 6,3,3 \rangle\) & 41 &   12-5.2 & \(\{(F_4,F_5,F_6,F_7,F_{10},F_{11}),(F_1,F_8,F_9),\) &\(F_3F_4,F_3F_6,F_4F_{12},F_6F_{12}\)\\
&&&&&&~~~~~~~~\(  (F_2,F_3,F_{12})  \}\)&  \\
12 & 5 &128& \(\langle 6,3,3 \rangle\) & 41 &   12-5.2 & \(\{(F_1,F_2,F_5,F_8,F_{10},F_{11}), (F_3,F_4,F_9),\) &\(F_3F_6,F_3F_{12},F_4F_6,F_4F_{12}\)\\
&&&&&&~~~~~~~~\( (F_6,F_7,F_{12})  \}\)&  \\
12 & 5 &128& \(\langle 5,5,2 \rangle\) & 41 &   12-5.2 & \(\{(F_1,F_3,F_5,F_8,F_{12}),(F_2,F_4,F_6,F_7,F_9),\) &\(F_3F_4,F_3F_6,F_4F_{12},F_6F_{12}\) \\
&&&&&&~~~~~~~~\( (F_{10},F_{11})  \}\)& \\
12 & 5 &128& \(\langle 7,3,2 \rangle\) & 37 &   12-5.3 & \(\{(F_2,F_3,F_4,F_6,F_9,F_{11},F_{12}), (F_1,F_7,F_8),\) &\(F_1F_2,F_1F_3,F_2F_8,F_3F_8\)\\
&&&&&&~~~~~~~~\( (F_5,F_{10})  \}\)&\\
12 & 5 &128& \(\langle 8,2,2 \rangle\) & 36 &   12-5.1 & \(\{(F_1,F_6,F_7,F_8,F_9,F_{10},F_{11},F_{12}), (F_2,F_3),\) &none\\
&&&&&&~~~~~~~~\( (F_4,F_5)  \}\)&\\
\hline
\end{tabular}
}
   
\vspace{0.1cm}

\(^a\)number of estimable interactions

\(^b\)interactions with constituent factors from different factor groups which are aliased with another interaction
   
\end{minipage}
\end{table}

\section{Concluding Remarks}
\label{disc}

This work develops a practitioner led approach to the construction of blocked \(2^n\) and \(2^{n-p}\) designs to estimate main effects and selected two-factor interactions. The generator matrix provides immediate and illuminating information on the capacity of a design to estimate the chosen effects. Hence,  the construction focus on the generator matrix, rather than on the block defining contrast sub-group, facilitates the building of bespoke designs. The method is especially advantageous and convenient when the block size is relatively small. 
It is notable that the most useful blocked \(2^n\) factorials are not necessarily those which maximise the number of estimable interactions. Likewise, the most useful blocked \(2^{n-p}\) factorials are not necessarily based on minimum aberration \(2^{n-p}\) fractions.

A novel aspect of the work is the use of concepts from graph theory. The representation of a blocked \(2^n\) factorial by a requirements graph with proper vertex colouring enables exploitation of results on chromatic numbers and proper vertex colourings.  

The design class with blocks of size four is an important one and is the class which gives the most construction challenges.
The difficulty is, in part, due to the large number of factorial effects that are confounded with blocks. The other feature of blocks of size four that contributes to the design limitations is the existence of small configurations of interactions that cannot be accommodated: any interaction set containing all six pairwise interactions between four factors will correspond to a requirements graph with chromatic number at least four. No design can be found to cover all main effects and such an interaction set. For fractional designs The tables of designs provided for designs in blocks of size four  should form a useful resource.

\appendix

\section{Appendix}\label{app}

The \(2^{n-p}\) factorials listed in Table \ref{tab:fg} are obtained from Chen {\it et al.} (1993) and Block and Mee (2005).


\begin{table}
\centering
\begin{minipage}{11cm}
\caption{\(2^{n-p}\) Factorials} \label{tab:fg}

\begin{tabular}{cl}\hline
fractional & defining words   \\
factorial & \\ \hline
5-1.1 &\(F_1F_2F_3F_4F_5\)\\
6-1.1 &\(F_1F_2F_3F_4F_5F_6\)\\
6-1.2 &\(F_1F_2F_3F_4F_6\)\\
7-1.1 &\(F_1F_2F_3F_4F_5F_6F_7\)\\
7-1.2 &\(F_1F_2F_3F_4F_5F_7\)\\
7-2.1 & \(F_1F_2F_3F_6,F_1F_2F_4F_5F_7,\)\\
8-1.1 &\(F_1F_2F_3F_4F_5F_6F_7F_8\)\\
8-1.2 &\(F_1F_2F_3F_4F_5F_6F_8\)\\
8-2.1 & \(F_1F_2F_3F_4F_7, F_1F_2F_5F_6F_8\)\\
9-2.1 &\(F_1F_2F_3F_4F_5F_8, F_1F_2F_3F_6F_7F_9\)\\
9-2.2 &\(F_1F_2F_3F_4F_8, F_1F_2F_5F_6F_7F_9\)\\
9-3.1 & \(F_1F_2F_3F_7,F_1F_2F_4F_5F_8,F_1F_3F_4F_6F_9\)\\
10-3.1 &\(F_1F_2F_3F_4F_5F_8, F_1F_2F_3F_6F_7F_9,F_1F_2F_4F_6F_{10}\)\\
11-4.1 &\(F_1F_2F_3F_4F_5F_8,F_1F_2F_3F_6F_7F_9, F_1F_2F_4F_6F_{10}, F_1F_3F_5F_7F_{11}\)\\
12-5.1 &\(F_1F_2F_3F_4F_5F_8,F_1F_2F_3F_6F_7F_9, F_1F_2F_4F_6F_{10},F_1F_3F_5F_7F_{11},F_1F_4F_5F_6F_7F_{12}\)\\
12-5.2 &\(F_1F_2F_3F_4F_5F_8,F_1F_2F_3F_6F_7F_9, F_1F_2F_4F_6F_{10},F_1F_3F_5F_7F_{11},F_3F_4F_6F_{12}\)\\
12-5.3 & \(F_1F_2F_3F_8,F_1F_2F_4F_5F_9, F_1F_3F_4F_6F_{10},F_2F_3F_5F_7F_{11},F_4F_5F_6F_7F_{12}\)\\
\hline
\end{tabular}

\end{minipage}

\end{table}

{}
\end{document}